\input amstex
\documentstyle{amams} % input Annals of Mathematics macros.
\input amssym.def
\input amssym.tex

\document
\annalsline{153}{2001}
\received{September 3, 1999}
 \startingpage{259}
\font\tenrm=cmr10

\catcode`\@=11
\font\twelvemsb=msbm10 scaled 1100

%\font\ninemsb=msbm7 scaled 1100%msbm9
\font\ninemsb=msbm10 scaled 800
\newfam\msbfam
\textfont\msbfam=\twelvemsb  \scriptfont\msbfam=\ninemsb
  \scriptscriptfont\msbfam=\ninemsb
\def\msb@{\hexnumber@\msbfam}
\def\Bbb{\relax\ifmmode\let\next\Bbb@\else
 \def\next{\errmessage{Use \string\Bbb\space only in math
mode}}\fi\next}
\def\Bbb@#1{{\Bbb@@{#1}}}
\def\Bbb@@#1{\fam\msbfam#1}
\catcode`\@=12

 \catcode`\@=11
\font\twelveeuf=eufm10 scaled 1100
\font\teneuf=eufm10
\font\nineeuf=eufm7 scaled 1100%eufm9
\newfam\euffam
\textfont\euffam=\twelveeuf  \scriptfont\euffam=\teneuf
  \scriptscriptfont\euffam=\nineeuf
\def\euf@{\hexnumber@\euffam}
\def\frak{\relax\ifmmode\let\next\frak@\else
 \def\next{\errmessage{Use \string\frak\space only in math
mode}}\fi\next}
\def\frak@#1{{\frak@@{#1}}}
\def\frak@@#1{\fam\euffam#1}
\catcode`\@=12

%--------------- Author macros ---------------

 \def\famzero{\rm }
\def\arccos{\mathop{\famzero arccos}\nolimits}
\def\arcsin{\mathop{\famzero arcsin}\nolimits}

\def\arg{\mathop{\famzero arg}\nolimits}

\def\det{\mathop{\famzero det}}

\def\exp{\mathop{\famzero exp}\nolimits}

\def\log{\mathop{\famzero log}\nolimits}
\def\max{\mathop{\famzero max}}
\def\min{\mathop{\famzero min}}

\def\sin{\mathop{\famzero sin}\nolimits}

\define\pperm{\Bbb P_{\text{perm},N}}
\define\pplanchN{\Bbb P_{\text{Plan},N}}
\define\pplancha{\Bbb P_{\text{Plan}}^\alpha}
\define\pcharN{\Bbb P_{\text{Ch,M,N}}}
\define\pchara{\Bbb P_{\text{Ch},M}^\alpha}
\define\Ai{\text{Ai}\,}
\define\KMe{K_{\,\text{Me},M}^{K,q}}
\define\KCh{K_{\,\text{Ch},M}^\alpha}
\define\eplana{\Bbb E_{\,\text{Plan}}^\alpha}
\define\emeix{\Bbb E_{\,\text{Me},N,N}^{\alpha/N^2}}
\define\echara{\Bbb E_{\,\text{Ch},M}^\alpha}
\define\echarN{\Bbb E_{\,\text{Ch},M,N}}
\define\pwMN{\Bbb P_{\text{W},M,N}}
\define\pwMa{\Bbb P_{\text{W},M}^\alpha}

%-------------- Author entries --------------------
\title{Discrete orthogonal polynomial ensembles\\ and the Plancherel
measure} %Article title
\shorttitle{Ensembles  and the Plancherel
measure} 
% Acknowledgements: Please enter all acknowledgements here.
\acknowledgements{}
 \author{Kurt Johansson}
   \institutions{Royal Institute of Technology, Stockholm, Sweden\\
 {\eightpoint {\it E-mail address\/}: kurtj\@math.kth.se}}

\bigbreak \centerline{\bf Abstract} 
\bigbreak
We consider discrete orthogonal polynomial ensembles
which are discrete analogues of the orthogonal polynomial ensembles in
random matrix theory. These ensembles occur in certain problems in 
combinatorial probability and can be thought of as probability measures
on partitions. The Meixner ensemble is related to a two-dimensional
directed growth model, and the Charlier ensemble is related to the lengths of 
weakly increasing subsequences in random words. The Krawtchouk ensemble
occurs in connection with zig-zag paths in random domino tilings of the
Aztec diamond, and also in a certain simplified directed first-passage
percolation model.
We use the Charlier ensemble 
to investigate the asymptotics of weakly increasing 
subsequences in random words and to prove a conjecture of Tracy and Widom.
As a limit of the Meixner ensemble or the Charlier ensemble we obtain the
Plancherel measure on partitions, and using this we prove a conjecture of
Baik, Deift and Johansson that under the Plancherel measure, the distribution
of the lengths of the first $k$ rows in the partition, appropriately scaled,
converges to the asymptotic joint distribution for the $k$ largest eigenvalues
of a random matrix from the Gaussian Unitary Ensemble. In this problem a 
certain discrete kernel, which we call the discrete Bessel kernel, plays
an important role.

\section{Introduction and results}

During the last years there has been a lot of activity around the
problem of the distribution of the length of a longest increasing
subsequence of a random permutation, its generalizations and their
connection with random matrices, see for example [Ge], [Ra], 
[BDJ1], [Jo3], [Ok], [BR2], [Bi],  
and also [AD] for connections with patience and
the history of the problem. Let $\pi$ be a random permutation from the
symmetric group $S_N$ with uniform distribution $\pperm$ and let
$L(\pi)$ denote the length of a longest increasing subsequence in
$\pi$.
It is proved by Baik, Deift and Johansson in [BDJ1] that
$$
\lim_{N\to\infty}\pperm[L(\pi)\le 2\sqrt{N}+tN^{1/6}]=F(t),\tag 1.1
$$
where $F(t)$ is the {\it Tracy-Widom distribution}, (1.5)  
for the appropriately
scaled largest eigenvalue of a random $M\times M$ matrix from the
{\it Gaussian Unitary Ensemble} (GUE) in the limit $M\to\infty$, see
[TW1]. The probability
density function on $\Bbb R^M$ for the $M$ eigenvalues $x_1,\dots,x_M$
of an $M\times M$ GUE matrix is
$$
\phi_{{\rm GUE},M}(x)=\frac 1{Z_M}\prod_{1\le i<j\le
M}(x_i-x_j)^2\prod_{j=1}^M e^{-x_j^2},\tag 1.2
$$  
where $Z_M=(2\pi)^{M/2}2^{-M^2/2}\prod_{j=1}^M(j!)^{-1}$. This
probability density can be analyzed using the Hermite polynomials,
which are orthogonal with respect to the weight $\exp(-x^2)$ occurring
in (1.2). Using standard techniques from random matrix theory, see
[Me] or [TW2], we can write
$$
\Bbb P_{{\rm GUE},M}\left[\max_{1\le k\le M} x_k\le
\sqrt{2M}+\frac{t}{\sqrt{2}M^{1/6}}\right]=\det(I-\Cal K_M)\bigr|_{L^2(t,\infty)},\tag
1.3
$$
where
$$
\Cal K_M(\xi,\eta)=\frac{1}{\sqrt{2}M^{1/6}} 
K_M\left(\sqrt{2M}+\frac{\xi}{\sqrt{2}M^{1/6}}, 
\sqrt{2M}+\frac{\eta}{\sqrt{2}M^{1/6}}\right).
$$
Here $K_M$ is the {\it Hermite kernel},
$$
K_M(x,y)=\frac{\kappa_{M-1}}{\kappa_M}\frac{h_M(x)h_{M-1}(y)-
h_{M-1}(x)h_M(y)}{x-y}e^{-(x^2+y^2)/2}
$$
with $h_m(x)=\kappa_m x^m+\dots$,
$\int_{\Bbb R}h_n(x)h_m(x)\exp(-x^2)dx=\delta_{nm}$,
the normalized Hermite polynomials. It follows from
standard asymptotic results for Hermite polynomials that
$$
\lim_{M\to\infty}\Cal K_M(\xi,\eta)=A(\xi,\eta)\doteq
\frac{\Ai(\xi)\Ai'(\eta)-\Ai'(\xi)\Ai(\eta)}{\xi-\eta},\tag 1.4
$$
the {\it Airy kernel}, and also that the Fredholm determinant in the
right-hand side of (1.3) converges to
$$
F(t)=\det(I-A)\bigr|_{L^2(t,\infty)}=\sum_{k=0}^\infty\frac{(-1)^k}{k!}
\int_{(t,\infty)^k}\det[A(\xi_i,\xi_j)]_{i,j=1}^kd^k\xi,\tag 1.5
$$
the Tracy-Widom distribution.

The problem of the length of the longest increasing subsequence in a
random permutation is closely related to the so called Plancherel
measure on partitions, which occurs as a natural probability measure
on the set of all equivalence classes of irreducible representations
of the symmetric group. Let
$\lambda=(\lambda_1,\lambda_2,\dots,\lambda_{\ell},0,0,\dots),
\lambda_1\ge\lambda_2\ge\dots\ge\lambda_{\ell}\ge 1$,
$\sum_j\lambda_j=N$, be a partition of $N$, which can be represented
in the usual way by a Young diagram with $\ell$ rows and $\lambda_j$
boxes in the $j^{\rm th}$ row, see e.g. [Sa], [Fu]. Let $f^\lambda$ be the
number of standard Young tableaux of shape $\lambda$. The {\it Plancherel
measure} assigns to $\lambda$ the probability
$$
\pplanchN[\{\lambda\}]=\frac{(f^\lambda)^2}{N!}.\tag 1.6
$$
The probability measure (1.6) is the push-forward of the uniform
distribution on $S_N$ by 
the Robinson-Schensted-Knuth (RSK)-correspondence, see
e.g. [Sa] or [Fu], which maps a permutation $\pi$ to a pair of
standard Young tableaux of the same shape $\lambda$, and 
the length $\lambda_1$ of the first row 
is equal to $L(\pi)$. Thus,
the length of the first row behaves in the limit as $N\to\infty$, as
the largest eigenvalue of a GUE matrix. It was proved in [BDJ2] that
the distribution of the rescaled length of the second row, $\pplanchN
[\lambda_2\le 2\sqrt{N}+tN^{1/6}]$, converges to the Tracy-Widom
distribution for the second largest eigenvalue of a GUE matrix, [TW2],
and it was conjectured that the analogous result holds for the $k^{\rm th}$
row. This conjecture will be proved in the present paper. It has
recently been independently proved by Borodin, Okounkov and Olshanski,
[BOO], see below. The conjecture also follows from the result by
Okounkov in [Ok]. His proof uses interesting geometric/combinatorial
methods. There are many earlier indications of connections between the
Plancherel measure and random matrices for instance in the work of
Regev, [Re], and Kerov, [Ke1], [Ke2].

Another measure on partitions, coming from pairs of semi-standard\break
tableaux, arises in [Jo3], where a certain random growth model is
investigated. This measure relates to a discrete Coulomb gas on $\Bbb
N$ of the form
$$
\frac 1{Z_M}\prod_{1\le i<j\le M}(h_i-h_j)^2\prod_{j=1}^Mw(h_j),\quad
h\in \Bbb N^M,\tag 1.7
$$
where $Z_M$ is a normalization constant. The weight
$w(x)=\binom{x+K-1}{x}q^x$, is the weight function on $\Bbb N$
for the Meixner polynomials, $m_n^{K,q}(x)$, see [NSU]. This measure
on $\Bbb N^M$ can be analyzed using the {\it Meixner kernel}
$$\multline
\KMe
(x,y)\\ =\frac{-q}{(1-q)d_{M-1}^2}\frac{m_M(x)m_{M-1}(y)-m_{M-1}(x)m_M(y)}
{x-y}(w(x)w(y))^{1/2},\endmultline \tag 1.8
$$
with $d_n=n!(n+K-1)!(1-q)^{-K}q^{-n}[(K-1)!]^{-1}$, in much the same
way as (1.2) is analyzed using the Hermite kernel. The Meixner kernel
occurs in connection with probability measures on partitions also in
the work of Borodin and Olshanski, [BO1]. The connection between
certain measures on partitions and discrete Coulomb gases with their
associated orthogonal polynomials is central in the present paper, and
give them a very interesting statistical mechanical interpretation
very similar to Dyson's Coulomb gas picture of 
the eigenvalues of random matrices. The difference is that in (1.7) we
have a Coulomb gas on the integer lattice instead of on the real line.
Other
statistical mechanical aspects of measures on partitions have been
investigated by Vershik, see [Ve] and references therein.
We will refer to (1.7) as a {\it discrete orthogonal
polynomial ensemble}. We will also be concerned with the cases
$w(x)=\alpha^xe^{-\alpha}/x!$, $x\in\Bbb N$, the {\it Charlier
ensemble},
$w(x)=\binom{N}{x}p^xq^{N-x}$, $x\in\{0,\dots,N\}$, the {\it 
Krawtchouk ensemble} and $w(x)$ given by (5.19), the {\it Hahn ensemble}.

Consider the {\it Poissonized Plancherel measure},
$$
\pplancha[\{\lambda\}]=e^{-\alpha}\sum_{N=0}^\infty\pplanchN[\{\lambda\}]
\frac {\alpha^N}{N!},\tag 1.9
$$
on the set of all partitions, $\pplanchN[\{\lambda\}]=0$ if
$\sum_j\lambda_j\neq N$. We will prove that this measure is a limit as
$q\to 0$ of the Meixner ensemble. The Meixner kernel (1.8) converges
in this limit, ($q=\alpha/M^2$, $K=1$, $M\to\infty$), 
to the {\it discrete Bessel kernel}
$$
B^\alpha(x,y)=\sqrt{\alpha}\frac{J_x(2\sqrt{\alpha})J_{y+1}(2\sqrt{\alpha})
-J_{x+1}(2\sqrt{\alpha})J_y(2\sqrt{\alpha})}{x-y}.\tag 1.10
$$
This result can be used to give a new proof of (1.1), and also to
verify the $k^{\rm th}$ row conjecture of [BDJ2], as well as to obtain
asymptotic results in the ``bulk'' of the Young diagram. These results
have recently been independently obtained by Borodin, Okounkov and
Olshanski, [BOO], as a limiting case of the results in [BO1]. See
the paper [BO2] for a discussion of 
the connections between [BOO] and the present
paper.

The results for the Poissonized Plancherel measure can also be
obtained as a limit of the Charlier ensemble. This ensemble arises in
the problem of the distribution of the length of a longest weakly
increasing subsequence in a random word which will be studied below.
The random word problem has recently been investigated by Tracy and
Widom, [TW3], using Toeplitz determinants and Painlev\'e equations,
see also [AD].

Before stating our results precisely we must introduce some
notation. Let 
$$
\Cal P=\{\lambda\in\Bbb N^{\Bbb Z_+}\,;\,\lambda_1\ge\lambda_2\ge\dots
\,\, \text{and}\,\,\sum_j\lambda_j<\infty\}
$$
denote the set of all partitions, and $\Cal P^{(N)}=\{\lambda\in\Cal
P\,;\,\sum_j\lambda_j =N\}$, $N\ge 0$, the set of all partitions of
$N$. Set $\ell(\lambda)=\max\{k\,;\,\lambda_k>0\}$, the {\it length} of
$\lambda$. We will consider functions on $\Cal P$ of the following
form. Let $f:\Bbb Z\to\Bbb C$ be a bounded function which satisfies
$f(n)=1$ if $n<0$. For a given $L\ge 0$ we define $g:\Cal P\to\Bbb C$
by
$$
g(\lambda)=\prod_{i=1}^\infty f(\lambda_i+L-i).\tag 1.11
$$
We say that $g$ is generated by $f$. Let $\Cal G_L$ denote the set of
all functions $g$ obtained in this way and write
$c(g)=||f||_\infty$. Let $\Cal P_M=\{\lambda\in\Cal
P\,;\,\ell(\lambda)\le M\}$ and $\Cal P_M^{(N)}=\Cal P_M\cap\Cal
P^{(N)}$. We also define, for $M\ge 1$, $N\ge 0$,
$$\align
\Omega_M&=\{\lambda\in\Bbb N^M\,;\,\lambda_1\ge\lambda_2\ge\dots
\ge\lambda_M\},\\ 
\Omega_M^{(N)}&=\{\lambda\in\Omega_M\,;\,\sum_{j=1}^M\lambda_j=N\}.
\endalign$$
Note that there is a natural bijection between $\Cal P_M$ and
$\Omega_M$ (and $\Cal P_M^{(N)}$ and $\Omega_M^{(N)}$). If $M\ge L$,
$g\in\Cal G_L$ and $\lambda\in\Cal P_M$, then
$$
g(\lambda)=\prod_{i=1}^Mf(\lambda_i+L-i),\tag 1.12
$$
since $f(n)=1$ if $n<0$, and we take (1.12) as our definition of $g$
on $\Omega_M$.

For $m\ge 1$ and $\lambda\in\Cal P$ we define
$$
V_m(\lambda)=\prod_{1\le i<j\le m}(\lambda_i-\lambda_j+j-i),
$$
and
$$
W_m(\lambda)=\prod_{i=1}^m\frac 1{(\lambda_i+m-i)!}.
$$
According to a formula of Frobenius, see e.g. [Sa] or [Fu], the quantity
$f^\lambda$ above can be expressed as
$$
f^\lambda=N!V_{\ell(\lambda)}(\lambda)W_{\ell(\lambda)}(\lambda).\tag
1.13
$$
Let $q\in (0,1)$ and $N\ge M$. We define the {\it Meixner ensemble} on
$\Omega_M$ by
\vglue-9pt
$$
\Bbb P_{\text{Me},M,N}^q[\{\lambda\}]=
(1-q)^{MN}\prod_{j=0}^{M-1}\frac{(N-M)!}{j!(N-M+j)!}V_M(\lambda)^2
\prod_{i=1}^M\binom{\lambda_i+N-i}{\lambda_i+M-i}q^{\lambda_i}.\tag
1.14
$$
Note that if we make the change of variables $h_i=\lambda_i+M-i$ this
gives us the discrete Coulomb gas (1.7) with the Meixner weight
$w(x)=\binom {x+K-1}{x}q^x$, where $K=N-M+1$. For more about the
Meixner ensemble and its probabilistic interpretations see [Jo3].
We can now state our first theorem.

\nonumproclaim{Theorem 1.1} For any $g\in\Cal G_L$, $L\ge 0${\rm ,} and
$\alpha>0$
we have that
$$
\eplana[g]=\lim_{N\to\infty}\emeix[g].\tag 1.15
$$
\endproclaim

Thus the Poissonized Plancherel measure can be obtained as a limit of
the Meixner ensemble. The theorem will be proved in Section 2.

Next, we define the {\it Charlier ensemble} on $\Omega_M$,
which can be obtained as a limit of the
Meixner ensemble, see (3.1). Given $\alpha>0$ we define
$$
\pchara[\{\lambda\}]=\left(\prod_{j=1}^{M-1}\frac 1{j!}\right)V_M(\lambda)^2
W_M(\lambda)\prod_{i=1}^M\left[\left(\frac{\alpha}{M}\right)^{\lambda_i}
e^{-\alpha/M}\right] \tag 1.16
$$
on $\Omega_M$. Again, the change of variables $h_i=\lambda_i+M-i$
gives a discrete Coulomb gas, (1.7). The Poissonized Plancherel
measure can also be obtained as a limit of the Charlier ensemble.

\nonumproclaim{Theorem 1.2} For any $g\in\Cal G_L$, $L\ge 0${\rm ,} and
$\alpha>0${\rm ,}
$$
\eplana[g]=\lim_{M\to\infty}\echara[g].\tag 1.17
$$
\endproclaim

The Charlier ensemble has a probabilistic interpretation in terms of
random words, see Proposition 1.5. Since the Meixner and Charlier
ensembles both correspond to discrete orthogonal polynomial ensembles
they can be analyzed in a way similar to that in which the Hermite
ensemble (GUE) is analyzed. This makes it possible to prove the
following theorem, compare with [BOO].

\nonumproclaim{Theorem 1.3} Let $g\in\Cal G_L$, $L\ge 0${\rm ,} be generated by
$f${\rm ,} see {\rm (1.11),} and write $\phi=f-1${\rm .} Then{\rm ,}
$$
\eplana[g]=\sum_{k=0}^\infty\frac 1{k!}\sum_{h\in\Bbb N^k} 
\prod_{j=1}^k\phi(h_j)\det[B^\alpha(h_i-L,h_j-L)]_{i,j=1}^k,\tag 1.18
$$
where $B^\alpha$ is the discrete Bessel kernel{\rm , (1.10).}  Note that the
right\/{\rm -}\/hand side is the Fredholm determinant of the operator on
$\ell^2(\Bbb N)$ with kernel  $B^\alpha(x-L,\break y-L)\phi(y).$ 
\endproclaim

The theorem will be proved in Section 3. 

As an example we can take $\phi(t)=-\chi_{(n,\infty)}(s)$ and $L=0$. This gives
$$
\pplancha[\lambda_1\le n]=\det(I-B^\alpha)\bigl|_{\ell^2(\{n,n+1,\dots\})}.
$$
By Gessel's formula the left-hand side is also a certain Toeplitz
determinant, see e.g. [BDJ1], and hence we get an interesting identity 
between a Toeplitz determinant and a certain Fredholm determinant on a
discrete space. This formula has recently been generalized by Borodin
and Okounkov, [BoOk]. 

By letting $\alpha$ go to
infinity we can use (1.18) combined with de-Poisson\-ization techniques
to prove asymptotic properties of the Plancherel measure. In particular
the next theorem generalizes the results of [BDJ1] and [BDJ2]. Note,
however, that we do not prove convergence of moments of the
appropriately rescaled random variables. In Section 3 we will prove

\nonumproclaim{Theorem 1.4} Let $x^{(j)}$ denote the $j^{\rm th}$ largest
eigenvalue among the eigenvalues $x_1,\dots ,x_M$ of a random $M\times
M$ matrix from {\rm GUE} with measure {\rm (1.2).} There is a distribution
function $F(t_1,\dots,t_k)$ on $\Bbb R^k${\rm ,} see {\rm (3.48),} such that
\vglue-9pt
$$
\lim_{M\to\infty}\Bbb P_{{\rm GUE},M}\left[x^{(j)}\le\sqrt{2M}
+\frac{t_j}{\sqrt{2}M^{1/6}},j=1,\dots,k\right]=F(t_1,\dots,t_k),\tag 1.19
$$
\medbreak \noindent 
for   $(t_1,\ldots ,t_k)\in\Bbb R^k${\rm ,} and
\vglue-9pt
$$
\lim_{N\to\infty}\pplanchN[\lambda_j\le 2\sqrt{N}+t_jN^{1/6},j=1,\dots,k]=
F(t_1,\dots,t_k),\tag 1.20
$$
\medbreak \noindent 
for   $(t_1,\ldots , t_k)\in\Bbb R^k${\rm .}
\endproclaim

\vglue12pt 

We turn now to the random word problem. By a {\it word} of length $N$ on $M$
letters, $M,N\ge 1$, we mean a map
$w:\{1,\dots,N\}\to\{1,\dots,M\}$. Let $W_{M,N}$ denote the set of all
such words, and let $\pwMN[\cdot]$ be the uniform probability
distribution on $W_{M,N}$ where all $M^N$ words have the same
probability. A weakly increasing subsequence of $w$ is a subsequence
$w(i_1),\dots,w(i_m)$ such that $i_1<\dots<i_m$ and $w(i_1)\le\dots\le
w(i_m)$. Let $L(w)$ be the length of a longest weakly increasing subsequence
in $w$. The RSK-correspondence defines a bijection from $W_{M,N}$ to
the set of all pairs of Young tableaux $(P,Q)$ of the same shape
$\lambda\in\Cal P^{(N)}$, where $P$ is semistandard with elements in 
$\{1,\dots,M\}$ and $Q$ is standard with elements in
$\{1,\dots,N\}$. Under this correspondence $L(w)=\lambda_1$, the
length of the first row. Note that we must have $\ell(\lambda)\le M$,
so $\ell\in\Cal P_M^{(N)}$ which we can identify with
$\Omega_M^{(N)}$. In this way we get a map
$S:W_{M,N}\to\Omega_M^{(N)}$.

\vglue12pt 

\nonumproclaim{Proposition 1.5}\hskip-6pt The push\/{\rm -}\/forward of the uniform distribution
on $W_{M,N}$ by the map $S:W_{M,N}\to\Omega_M^{(N)}$ is
$$
\pwMN[S^{-1}(\lambda)]=\pcharN[\{\lambda\}]
\doteq\frac{N!}{M^N}\bigl(\prod_{j=1}^{M-1} \frac
1{j!}\bigr)V_M(\lambda)^2W_M(\lambda)\tag 1.21
$$
\medbreak \noindent on $\Omega_M^{(N)}${\rm .} The Poissonization of this measure is the
Charlier ensemble {\rm (1.16).} Consequently{\rm ,}
\vglue-9pt
$$
\pwMN[L(w)\le t]=\pcharN[\lambda_1\le t],\tag 1.22
$$
\medbreak \noindent and for the Poissonized word problem{\rm ,}
\vglue-9pt
$$
\pwMa[L(w)\le t]\doteq\sum_{N=0}^\infty e^{-\alpha}\frac{\alpha^N}{N!} 
\pwMN[L(w)\le t]=\pchara[\lambda_1\le t].\tag 1.23
$$
\endproclaim

\demo{Proof} See Section 4.\enddemo

\vglue8pt

The probability (1.23) can also be expressed as a Toeplitz determinant
using Gessel's formula, [Ge], see also [TW3]
and [BR1]. The formula (1.21) can be used to prove a conjecture by
Tracy and Widom, [TW3]. This conjecture says that the Poissonized
measure on $\Omega_M$ induced by the uniform distribution on words
converges, after appropriate rescaling, to the $M\times M$ GUE measure
(1.2). In Section 4 we will prove

\vglue8pt 

\nonumproclaim{Theorem 1.6} Let $g$ be a continuous function on $\Bbb
R^M${\rm .} Then 
\vglue-9pt 
$$\multline \lim_{N\to\infty}\echarN\biggl[g\biggl(
\frac{\lambda_1-N/M}{\sqrt{2N/M}},\dots,
\frac{\lambda_M-N/M}{\sqrt{2N/M}}\biggr)\biggr]\\ =M!\sqrt{\pi M}\int
_{\Bbb A_M} g(x)\phi_{{\rm GUE},M}(x)dx_1\dots dx_{M-1},
\endmultline\tag
1.24$$
\smallbreak\noindent where $\Bbb A_M=\{x\in\Bbb R^M\,;\,\text{$x_1>\dots>x_M$
 and $x_1+\dots+x_M=0$}\}${\rm .} Furthermore\vglue-9pt
$$ 
\multline
 \lim_{\alpha\to\infty}\echara\biggl[g\biggl(
\frac{\lambda_1-\alpha/M}{\sqrt{2\alpha/M}},\dots,
\frac{\lambda_M-\alpha/M}{\sqrt{2\alpha/M}}\biggr)\biggr] \\ =
M!\int_{\{x\in\Bbb R^M\,;\,x_1>\dots>x_M\}}g(x)\phi_{{\rm GUE},M}(x)
d^Mx.
\endmultline\tag1.25$$
\endproclaim  

 \vglue8pt 

The case when $g$ only depends on $\lambda_1$ has been proved in [TW3]
using very different methods. 

The formula (1.23) can be used to
analyze the asymptotics of the random variable $L(w)$ on $W_{M,N}$ as
both $M$ and $N$ go to infinity.

\vglue8pt 

\nonumproclaim{Theorem 1.7} Let $F(t)$ be the Tracy\/{\rm -}\/Widom distribution
function {\rm (1.5).} Then{\rm ,} for all $t\in\Bbb R${\rm ,}
$$
\lim_{\alpha\to\infty}\pwMa\left[L(w)\le\frac{\alpha}M+2\sqrt{\alpha}+\biggl( 
1+\frac{\sqrt{\alpha}}{M}\biggr)^{2/3}\alpha^{1/6}t\right]=F(t).\tag 1.26
$$
Assume that $M=M(N)\to\infty$ as $N\to \infty$ in such a way that
$(\log N)^{1/6}/M(N)\break\to 0${\rm .} Then{\rm ,} for all $t\in\Bbb R${\rm ,}
$$ 
\lim_{N\to\infty}\pwMN\left[L(w)\le\frac{N}M+2\sqrt{N}+\biggl( 
1+\frac{\sqrt{N}}{M}\biggr)^{2/3}N^{1/6}t\right]=F(t).\tag 1.27
$$
\endproclaim

\vglue8pt 

\demo{Proof} See Section 4. \enddemo

Note that when $M\gg\alpha$, the leading order of the mean goes like
$2\sqrt{\alpha}$ and the standard deviation like $\alpha^{1/6}$ just
as for random permutations. When\break $M\ll\alpha$, we expect from (1.3)
and (1.25) that
$$\align L(w)&=\lambda_1\approx\alpha/M+\sqrt{2\alpha/M}(\sqrt{2M}
+t/\sqrt{2}M^{1/6} )\\
&=\alpha/M+2\sqrt{\alpha}+t\sqrt{\alpha}/M^{2/3},\\
\noalign{\noindent which fits perfectly with (1.26).}\endalign$$

In Section 5 we will consider two problems in combinatorial
probability that relate to the Krawtchouk ensemble, namely Sepp\"al\"ainen's
simplified model of directed first-passage percolation and zig-zag paths
in random domino tilings of the Aztec diamond introduced by Elkies,
Kuperberg, Larsen and Propp. Since both problems require some definitions we
will not state the results here. A third problem, random tilings of a
hexagon by rhombi, which is related to the Hahn ensemble will also be
discussed briefly.

\section{The Plancherel measure as a limit of the Meixner ensemble }

The setting is the same as in [Jo3]. Let $\Cal M_N$ denote the set of
all $N\times N$ matrices with elements in $\Bbb N$. We define a
probability measure, $\Bbb P_N^q[\cdot]$
on $\Cal M_N$ by letting each element $a_{ij}$ in
$A\in\Cal M_N$ be geometrically distributed with parameter $q\in
(0,1)$, and requiring all elements to be independent. Then
$$
\Bbb P_N^q[A]=(1-q)^{N^2}q^{\Sigma(A)},\tag 2.1
$$
$A\in\Cal M_N$, where $\Sigma(A)=\sum_{i,j=1}^N a_{ij}$. Let $\Cal
M_N(k)$ denote the set of all $A$ in $\Cal M_N$ for which
$\Sigma(A)=k$. Note that by (2.1) all matrices in $\Cal
M_N(k)$ have the same probability. Furthermore we let $\tilde{\Cal
M}_N(k)$ be the set of all matrices $A$ in $\Cal M_N(k)$ for which
$\sum_i a_{ij}\le 1$ for each $j$ and 
$\sum_j a_{ij}\le 1$ for each $i$; $\tilde{\Cal M}_N=\cup_k \tilde{\Cal
M}_N(k)$. By taking the appropriate submatrix of $A\in \tilde{\Cal
M}_N(k)$ we get a permutation matrix and hence a unique
permutation. This defines a map $R:\tilde{\Cal M}_N(k)\to S_k$, where
$S_k$ is the $k^{\rm th}$ symmetric group. Note that if $q$ is very small a
typical element in $\Cal M_N$ belongs to $\tilde{\Cal M}_N(k)$ for
some $k$. This is the crucial observation for what follows. The
RSK-correspondence defines a map $K:\Cal M_N(k)\to\Cal P^{(k)}$, and
also a map $S:S_k\to\Cal P^{(k)}$. The number of elements in $S_k$
that are mapped to the same $\lambda$ equals $(f^\lambda)^2$. It is
not difficult to see that if $A\in \tilde{\Cal M}_N(k)$ then
$K(A)=S(R(A))$. Let $g\in\Cal G_L$. It is proved in [Jo3] that
$$
\Bbb E_N^q[g(K(A))]=\Bbb E_{\text{Me},N,N}^q[g(\lambda)].\tag
2.2 
$$

With these preparations we are ready for the
\demo{Proof of Theorem {\rm 1.1}} By (2.2) we see that in order to prove
(1.15) it suffices to show that \vglue-9pt
$$
\lim_{N\to\infty}\Bbb E_N^{\alpha/N^2}[g(K(A))]=\eplana[g].\tag 2.3
$$
\smallbreak \noindent Note that $\Bbb P_N^q[\{A\}\bigr|\Sigma(A)=k]=1/\#\Cal M_N(k)$, where
$\#\Cal M_N(k)=\binom{N^2-1+m}{m}$,
and
$\Bbb P_N^q[\Sigma(A)=k]=\#\Cal M_N(k)(1-q)^{N^2}q^k$, by
(2.1). Thus
$$\align
&\Bbb E_N^q[g(K(A))\chi_{\tilde{\Cal M}_N}(A)]=
\sum_{k=0}^\infty\Bbb E_N^q[g(K(A))\chi_{\tilde{\Cal M}_N}(A)\bigr|\Sigma(A)=k]
\Bbb P_N^q[\Sigma(A)=k]\tag 2.4\\ \noalign{\vskip5pt}
&\qquad =\sum_{k=0}^\infty\sum_{A\in\tilde M_N(k)}g(K(A))(1-q)^{N^2}q^k\\
\noalign{\vskip5pt}
&\qquad =(1-q)^{N^2}\sum_{k=0}^\infty q^k\sum_{\lambda\in\Cal P^{(k)}}
g(\lambda)\#\{A\in\tilde M_N(k)\,;\,K(A)=\lambda\}. \\ \noalign{\vskip5pt}
\endalign
$$
The number of matrices in $\tilde M_N(k)$ which are mapped to the same
permutation by $R$ is $\binom{N}{k}^2$, since there are $\binom{N}{k}$
ways of choosing the rows and $\binom{N}{k}$ ways of choosing the
columns that select the submatrix. Since $K=S\circ R$ we obtain
$$
\#\{A\in\tilde M_N(k)\,;\,K(A)=\lambda\}=\binom{N}{k}^2(f^\lambda)^2.
$$
Together with (2.4) this yields
$$\align
\Bbb E_N^q[g(K(A))\chi_{\tilde{\Cal M}_N}(A)]&=
(1-q)^{N^2}\sum_{k=0}^\infty \frac{q^k}{k!}\frac{N!^2}{(N-k)!^2}
\sum_{\lambda\in\Cal P^{(k)}}
g(\lambda)\frac{(f^\lambda)^2}{k!}\\ \noalign{\vskip5pt}
&=(1-q)^{N^2}\sum_{k=0}^\infty \frac{q^k}{k!}\frac{N!^2}{(N-k)!^2}
\Bbb E_{\text{Plan},k}[g]\\ \noalign{\vskip5pt}
&=(1-\alpha/N^2)^{N^2}\sum_{k=0}^\infty\frac{\alpha^k}{k!}\biggl(
\frac{N!}{N^k(N-k)!}\biggr)^2\Bbb E_{\text{Plan},k}[g],\\ \noalign{\vskip5pt}
\endalign
$$
if we pick $q=\alpha/N^2$. Since $N!(N^k(N-k)!)^{-1}\le 1$ and
converges to 1 as $N\to\infty$ for each fixed $k$
and furthermore $\Bbb E_{\text{Plan},k}[g]\le c(g)^{\max(L,k)}$, 
it follows from the dominated convergence theorem that
$$
\lim_{N\to\infty}\Bbb E_N^{\alpha/N^2}[g(K(A))\chi_{\tilde{\Cal M}_N}(A)]=
\Bbb E_{\;\text{Plan}}^\alpha[g].\tag 2.5
$$
To deduce (2.3) from (2.5) we have to show that if $\tilde{\Cal
M}_N^\ast= \Cal M_N\setminus\tilde{\Cal M}_N$, then
$$
\Bbb E_N^{\alpha/N^2}[g(K(A))\chi_ {\tilde{\Cal
M}_N^\ast}(A)]=0.\tag 2.6
$$
By the Cauchy-Schwarz' inequality, the left-hand side of (2.6) is
$$
\le\Bbb E_N^{\alpha/N^2}[g(K(A))^2]^{1/2}
\Bbb P_N^{\alpha/N^2}[\tilde{\Cal
M}_N^\ast]^{1/2}.\tag 2.7
$$
If $\lambda=K(A)$, then $\ell(\lambda)\le\Sigma(A)$ and from the
definition (1.11) of $g$ it follows that
$$
|g(K(A))|\le c(g)^{\max(L,\ell(\lambda))}\le c(g)^{L+\Sigma(A)}.
$$
Thus,
$$
\Bbb E_N^{\alpha/N^2}[g(K(A))^2]\le c(g)^{2L}\sum_{k=0}^\infty
c(g)^{2k}\Bbb P_N^{\alpha/N^2}[\Sigma(A)=k].
$$
Since,
$$
\Bbb P_N^{\alpha/N^2}[\Sigma(A)=k]=\binom{N^2-1+k}{k}
\biggl(1-\frac{\alpha}{N^2}\biggr)^{N^2}\biggl(\frac{\alpha}{N^2}\biggr)^k\to
e^{-\alpha}\frac{\alpha^k}{k!}
$$
as $N\to\infty$, it is not hard to show that
$$
\Bbb E_N^{\alpha/N^2}[g(K(A))^2]\le C(\alpha,g),\tag 2.8
$$
for all $N\ge 1$, where $C(\alpha,g)$ depends only on $\alpha$ and
$c(g)$.

Next, we note that 
$$
\tilde{\Cal M}_N^\ast\subseteq\bigcup_{i=1}^N\left\{\sum_j a_{ij}\ge
2\right\}\cup \bigcup_{j=1}^N\left\{\sum_i a_{ij}\ge 2\right\}
$$
and hence
$$
\Bbb P_N^{\alpha/N^2}[\tilde{\Cal
M}_N^\ast]\le 2N\Bbb P_N^{\alpha/N^2}\left[\sum_j a_{ij}\ge 2\right].
$$
Since, $\Bbb P_N^{q}[\sum_j a_{ij}\ge
2]=1-(1-q)^N-N(1-q)^{N-1}$,
we obtain
$$
\Bbb P_N^{\alpha/N^2}[\tilde{\Cal
M}_N^\ast]\le\frac{C\alpha^2}{N}.
$$
Together with (2.7) and (2.8) this implies (2.6) and we are done.
\enddemo

\phantom{time to come home}

It is also possible to give a more direct proof based on the explicit
formulas similarly to what will be done with the Charlier ensemble in
the next section. Above we have emphasized the probabilistic and
geometric picture.

\section{The Plancherel measure as a limit of the Charlier
ensemble}

3.1. {\it The limit of the Charlier ensemble}.
The Charlier ensemble is defined by (1.16). It can be obtained as a
limit of the Meixner ensemble (1.14) by taking $q=\alpha/MN$ and
letting $N\to\infty$ with $M$ fixed. In this limit
$$\align
 &  (1-q)^{MN}\prod_{j=0}^{M-1}\frac{(N-M)!}{j!(N-M+j)!}\prod_{i=1}^M
\binom{\lambda_i+N-i}{\lambda_i+M-i}q^{\lambda_i}\tag 3.1\\&\hskip1in \to
\left(\prod_{j=1}^{M-1}\frac 1{j!}\right)W_M(\lambda)\prod_{i=1}^M
\left[\left(\frac{\alpha}{M}\right)^{\lambda_i}e^{-\alpha/M}\right],
\endalign$$ 
so we obtain (1.16). In   light of Theorem 1.1 we see that it is
reasonable to expect that the Poissonized Plancherel measure should be
the limit of the Charlier ensemble as $M\to\infty$. The interpretation
of the Charlier ensemble in connection with random words, Proposition
1.5, also supports this since a random word in the limit $M\to\infty$
is like a permutation (no letter is used twice), see also [TW3]. We
will give an analytical proof of Theorem 1.2 that does not use the
RSK-correspondence. 
We start with the following simple but important
lemma.

\nonumproclaim{Lemma 3.1} If $M\ge\ell(\lambda)${\rm ,} then
$$
V_M(\lambda)W_M(\lambda)=V_{\ell(\lambda)}(\lambda)W_{\ell(\lambda)}(\lambda).
\tag 3.2
$$
\endproclaim

\demo{Proof} We may assume that $M>\ell(\lambda)$. Note that, by
definition, $\lambda_i=0$ if $i>\ell(\lambda)$. Hence,
$$\align
V_M(\lambda)&=V_{\ell(\lambda)}(\lambda)\prod_{i=1}^{\ell(\lambda)}
\prod_{j=\ell(\lambda)+1}^M(\lambda_i+j-i)\prod_{\ell(\lambda)<i<j\le
M} (j-i)\\
&=V_{\ell(\lambda)}(\lambda)\prod_{i=1}^{\ell(\lambda)}
\frac{(\lambda_i+M-i)!}{(\lambda_i+\ell(\lambda)-i)!}
\prod_{\ell(\lambda)<i<j\le
M} (j-i).
\endalign
$$
Thus in order to prove (3.2) we must show that 
$$
\prod_{\ell(\lambda)<i<j\le
M} (j-i)=
\prod_{i=\ell(\lambda)+1}^M (\lambda_i+M-i)!,
$$
but this is immediate since $\lambda_i=0$ if $i>\ell(\lambda)$.
\enddemo

\demo{Proof of Theorem {\rm 1.2}}
It follows from the definition (1.12) of $g(\lambda)$ that
$$
|g(\lambda)|\le c(g)^{\max(\ell(\lambda),L)}\le c(g)^{\ell (\lambda)+L}.
\tag 3.3
$$
Let $\pcharN$ be defined by (1.21). Then,
$$
\align
\echara[g]&=\sum_{\lambda\in\Omega_M}
g(\lambda)\left(\prod_{j=1}^{M-1}\frac 1{j!}
\right)V_M(\lambda)^2W_M(\lambda)\prod_{j=1}^M\left[\left(\frac{\alpha}M\right)^
{\lambda_j}e^{-\alpha/M}\right]\tag 3.4\\
&=\sum_{N=0}^\infty e^{-\alpha}\frac{\alpha^N}{N!}
\sum_{\lambda\in\Omega_M^{(N)}} g(\lambda)\pcharN[\{\lambda\}].
\endalign
$$
Thus, by (3.3) and the fact that 
$\ell(\lambda)\le N$ if $\lambda\in\Omega_M^{(N)}$,
$$
\biggl|\sum_{\lambda\in\Omega_M^{(N)}}g(\lambda)\pcharN[\{\lambda\}]\biggr
|\le c(g)^{L+N},\tag 3.5
$$
since $\pcharN$ is a probability measure on $\Omega_M^{(N)}$.
Given $\varepsilon>0$ we can choose $K$ so large that
$$
\biggl|\sum_{N=K+1}^\infty e^{-\alpha}\frac{\alpha^N}{N!}c(g)^{L+N}
\biggr|\le\varepsilon.\tag 3.6
$$
Consequently,
$$
\biggl|\echara[g]-\sum_{N=0}^Ke^{-\alpha}\frac{\alpha^N}{N!}
\sum_{\lambda\in\Omega_M^{(N)}}g(\lambda)\pcharN[\{\lambda\}]\biggr|\le
\varepsilon.\tag 3.7
$$
If $M\ge K\ge N\ge\ell(\lambda)$, $\lambda\in\Omega_M^{(N)}$, we can identify
$\Omega_M^{(N)}$ with $\Cal P^{(N)}$ and use (3.2) to write
$$\align
\noalign{\vskip8pt}
&\tag 3.8\\
\noalign{\vskip-32pt}
\qquad &\sum_{\lambda\in\Omega_M^{(N)}}g(\lambda)\pcharN[\{\lambda\}]\\
&\hskip.5in=\sum_{\lambda\in\Omega_M^{(N)}}g(\lambda)N!V_{\ell(\lambda)}(\lambda)^2
W_{\ell(\lambda)}(\lambda)^2\prod_{j=1}^M\frac{(\lambda_j+M-j)!}{M^{\lambda_j}}
\prod_{j=1}^{M-1}\frac 1{j!}\\
&\hskip.5in=\sum_{\lambda\in\Cal P^{(N)}}g(\lambda)\pplanchN[\{\lambda\}]
\prod_{j=1}^{\ell(\lambda)}\frac{(\lambda_j+M-j)!}
{M^{\lambda_j}(M-j)!},
\endalign$$
where the last equality is a straightforward computation using the fact that 
$\lambda_j=0$ if $j>\ell(\lambda)$. Now,
$$
\prod_{j=1}^{\ell(\lambda)}\frac{(\lambda_j+M-j)!}
{M^{\lambda_j}(M-j)!}
=
\prod_{j=1}^{\ell(\lambda)}\left(1-\frac{j-1}{M}\right)\dots \left(1-\frac{j-\ell(\lambda)}M\right)
$$
which goes to 1 as $M\to\infty$ for a fixed $\lambda$. Since the sum in 
(3.7) is, for a fixed $K$, a sum over finitely many $\lambda$, we obtain 
$$\multline
 \lim_{M\to\infty}
\sum_{N=0}^Ke^{-\alpha}\frac{\alpha^N}{N!}
\sum_{\lambda\in\Omega_M^{(N)}}g(\lambda)\pcharN[\{\lambda\}] 
\\  =
\sum_{N=0}^Ke^{-\alpha}\frac{\alpha^N}{N!}
\sum_{\lambda\in\Cal P^{(N)}}g(\lambda)\pplanchN[\{\lambda\}].
\endmultline\tag 3.9$$
Using (3.6) and the fact that $\pplanchN$ is a probability measure on $\Cal P
^{(N)}$, we obtain
$$
\biggl|\eplana[g]-\sum_{N=0}^Ke^{-\alpha}\frac{\alpha^N}{N!}
\sum_{\lambda\in\Omega_M^{(N)}}g(\lambda)\pplanchN[\{\lambda\}]\biggr|\le
\varepsilon.\tag 3.10
$$
The theorem now follows from (3.7), (3.9) and (3.10).
\enddemo

\demo{{\rm 3.2.} Coulomb gas interpretation of the Plancherel measure}
As $M\to\infty$ the number of particles in the Coulomb gas representation 
of the Charlier ensemble goes to infinity, so a Coulomb gas interpretation 
of the Plancherel measure is not immediate. We will now show that we 
can actually approximate $\pplancha$ by a Coulomb gas with $K$ particles, 
which gives a good approximation if $K$ is chosen large enough (depending on
$\alpha$).

Consider the Poissonization of the 
restriction of the Plancherel measure to $\Cal
P_M^{(N)}$,
$$
F_M^\alpha[g]=e^{-\alpha}\sum_{N=0}^\infty
\frac{\alpha^N}{N!} \sum_{\lambda\in P^{(N)}_M} g(\lambda) {(f^\lambda)^2\over N!}
$$
for $g\in\Cal G_L$.
If $M\ge L$ it follows from (1.12), (1.13) and
Lemma 3.1 that
$$
F_M^\alpha[g]=e^{-\alpha}\sum_{\lambda\in
\Omega_M}g(\lambda)V_M(\lambda)^2W_M(\lambda)^2
\prod_{i=1}^M\alpha^{\lambda_i}.
$$
When $M$ is large, we expect that $F_M^\alpha[g]$ and $\eplana[g]$
should be close. \enddemo

\nonumproclaim{Lemma 3.2} Assume that $g\in\Cal G_L$ and let $d>0$ be
given{\rm .} There is a numerical constant $C$ such that if
$M\ge\max(L,\alpha\exp(d+1))${\rm ,} then
$$
\bigr|\eplana[g]-F_M^\alpha[g]\bigr|\le C(c(g)e^{-d})^M.\tag 3.11
$$
\endproclaim

\demo{Proof} Set 
$$
R_{N,M}[g]=\eplana[g]-F_M^\alpha[g]=\sum_{N=0}^\infty e^{-\alpha}
\frac{\alpha^N}{N!}
\sum_{\lambda\in\Cal
P^{(N)}\setminus \Cal P_M^{(N)}}g(\lambda)\frac{(f^\lambda)^2}{N!}.
$$
If $N\le M$, then $R_{N,M}[g]=0$ since then
$\ell(\lambda)\le\sum_i\lambda_i=N\le M$, so $\Cal
P^{(N)}=\Cal P_M^{(N)}$. If $N>M\ge L$, then $|g(\lambda)|\le c(g)^N$
since $\lambda_i=0$ if $i>N$. Thus,
$$
\align
\bigr|\eplana[g]-F_M^\alpha[g]\bigr|&\le\sum_{N=M+1}^\infty e^{-\alpha}
\frac{\alpha^N}{N!} \bigr|R_{N,M}[g]\bigr|\\
&\le\sum_{N=M+1}^\infty e^{-\alpha}
\frac{\alpha^N}{N!}c(g)^N.
\endalign$$
This last sum is estimated as follows. By Stirling's formula there is
a numerical constant $C$ such that $\exp(-\alpha)
\alpha^N/N!\le C\exp(-\alpha f(N/\alpha))$, where $f(x)=x\log
x+1-x$. If $N/\alpha\ge \exp(d+1)$, then $f(N/\alpha)\ge dN/\alpha$,
and so $\exp(-\alpha)
\alpha^N/N!\le \exp(-dN)$. The lemma is proved.
\enddemo

Recall from the introduction that $\Cal P_M$ can be naturally
identified with $\Omega_M$. For $K<M$ we define
$$
\Omega_{M,K}=\{\lambda\in\Omega_M\,;\,\lambda_{K+1}=\dots=\lambda_M=0\},
$$
and $\Omega_{M,K}^\ast=\Omega_M\setminus\Omega_{M,K}$. If $1\le j\le
M-K$ we set
$$
\Omega_{M,K}^\ast(j)=\{\lambda\in
\Omega_{M,K}^\ast\,;\,\lambda_{M+1-j}>0\,\,\text{but}\,\,\lambda_i=0,
M+1-j<i\le M\},
$$
so that $\Omega_{M,K}^\ast=\cup_{j=1}^{M-K}\Omega_{M,K}^\ast(j)$. The
next lemma asserts that $\ell(\lambda)$ is not too large for typical
$\lambda$ that we will consider.

\nonumproclaim{Lemma 3.3} Let $g\in\Cal G_L$ be generated by $f${\rm .} Assume
that $f$ satisfies
$$
0\le f(x)\le C_0f(x-1)\tag 3.12
$$
for all $x\in\Bbb Z$ and some constant $C_0${\rm .} Then
$$
e^{-\alpha}\sum_{\lambda\in\Omega_{M,K}^\ast(j)}g(\lambda)
V_M(\lambda)^2W_M(\lambda)^2 \prod_{i=1}^M\alpha^{\lambda_i}\le\frac
{(C_0\alpha)^{M-j+1}}{(M-j+1)!^2}F_M^\alpha[g].\tag 3.13
$$
\endproclaim

\demo{Proof} 
It will be most convenient to use the discrete Coulomb gas
representation.
Set $x_j=\lambda_{M+1-j}+j-1$, $j=1,\dots, M$ and let $\Delta_M(x)=
\prod_{1\le i<j\le M} (x_j-x_i)$ be the Vandermonde determinant. Also,
set
$A=\{x\in\Bbb N^M\,;\,0\le x_1<\dots<x_M\}$
and
$A_j=\{x\in A\,;\,\text{$x_i<i$ for $i<j$ and $x_j\ge j$}\}$,
$j=1,\dots, M$. Note that $\lambda\in\Omega_{M,K}^\ast(j)$ translates
into $x\in A_j$. If $x\in A_j$, then $x_i=i-1$ for $i=1,\dots,j-1$ and
we have the first \it hole \rm in the particle configuration $x$ at
$j-1$. Now,
$$
\sum_{\lambda\in\Omega_{M,K}^\ast(j)}g(\lambda)
V_M(\lambda)^2W_M(\lambda)^2 \prod_{i=1}^M\alpha^{\lambda_i}
=\sum_{x\in
A_j}\Delta_M(x)^2\prod_{i=1}^M\frac {\alpha^{x_i}}{x_i!^2}f(x_i+K-M).
\tag 3.14
$$
We want to show that, with high probability, the first hole must be
fairly close to $M$. Define $T_j:A_j\to A$ by
$T_j(x)=(x_1,\dots,x_{j-1},x_j-1,\dots,x_M-1)=x'$.
Clearly, $T_j:A_j\to T_j(A_j)$ is a bijection. Write
$$
L_M^\alpha(x)=\Delta_M(x)^2\prod_{i=1}^M\frac {\alpha^{x_i}}{x_i!^2}f(x_i+K-M).
$$
For $x\in A_j$,
$$
\biggl(\frac{\Delta_M(x)}
{\Delta_M(x')}\biggr)^2\prod_{i=1}^M\frac{(x_i'!)^2}
{(x_i!)^2}\alpha^{x_i-x_i'}=
\alpha^{M-j+1}\biggl(\frac{\Delta_M(x)}
{\Delta_M(x')}\biggr)^2\prod_{i=j}^M\frac{1}{x_i^2}.
$$
Since
$$
\prod_{i=j}^M\frac{1}{x_i^2}\le
\prod_{i=j}^M\frac{1}{i^2}=\biggl(\frac{(j-1)!}{M!}\biggr)^2
$$
and
$$
\frac{\Delta_M(x)}{\Delta_M(x')}=\prod_{k=j}^M\frac{x_k}{x_k-(j-1)}\le
\prod_{k=j}^M\frac{k}{k-(j-1)}=\binom{M}{j-1}
$$
if $x\in A_j$, we obtain, using our assumption on $f$,
$$
L_M^\alpha(x)\le
\frac{(C_0\alpha)^{M-j+1}}{(M-j+1)!^2}L_M^\alpha(T_j(x)).
$$
Inserting this into (3.14) yields
$$\multline
 e^{-\alpha}\sum_{\lambda\in\Omega_{M,K}^\ast(j)}g(\lambda)
V_M(\lambda)^2W_M(\lambda)^2 \prod_{i=1}^M\alpha^{\lambda_i}
=e^{-\alpha}\sum_{x\in A_j}L_M^\alpha(x)\\ 
\qquad\le
\frac{(C_0\alpha)^{M-j+1}}{(M-j+1)!^2}
e^{-\alpha}\sum_{x\in
A_j}L_M^\alpha(T_j(x))
\le\frac{(C_0\alpha)^{M-j+1}}{(M-j+1)!^2}F_M^\alpha[g],
\endmultline
$$
and the lemma is proved.
\enddemo
 
\nonumproclaim{Lemma 3.4} Let $g\in\Cal G_L$ be generated by $f$ which
satisfies {\rm (3.12).} Assume that
$M>K\ge\max(L,e\sqrt{2C_0\alpha})${\rm .} Then{\rm ,}
$$
\bigr|F_M^\alpha[g]-F_K^\alpha[g]\bigr|\le 2\biggl(\frac{C_0\alpha
e^2}{(K+1)^2}\biggr)^{K+1}F_M^\alpha[g].\tag 3.15
$$
\endproclaim

\demo{Proof}
If $\lambda\in\Omega_{M,K}$ then $\ell(\lambda)\le K<M$ and hence by
Lemma 3.1, (1.11) and the fact that $\Omega_{M,K}$ and $\Omega_K$ can
be identified we obtain
$$
e^{-\alpha}\sum_{\lambda\in\Omega_{M,K}}g(\lambda)V_M(\lambda)^2W_M(\lambda)^2
\prod_{i=1}^M\alpha^{\lambda_i}=F_K^\alpha[g].
$$
The left-hand side of (3.15) is
$$\align
&\le\sum_{j=1}^{M-K}
e^{-\alpha}\sum_{\lambda\in\Omega_{M,K}^\ast(j)}
g(\lambda)V_M(\lambda)^2W_M(\lambda)^2
\prod_{i=1}^M\alpha^{\lambda_i}\\
&\le\sum_{j=1}^{M-K}\frac{(C_0\alpha)^{M-j+1}}{(M-j+1)!^2}F_M^\alpha[g]\\
&\le\sum_{j=K+1}^\infty (C_0\alpha)^j\frac 1{j!^2}F_M^\alpha[g]\le 2
\biggl(\frac{C_0\alpha e^2}{(K+1)^2}\biggr)^{K+1}F_M^\alpha[g]
\endalign
$$
by (3.13). This completes the proof of the lemma.
\enddemo

We can now demonstrate how the Plancherel measure can be approximated by a 
Coulomb gas, (compare with the discussion in the Appendix in [BDJ1]).

\nonumproclaim{Proposition 3.5} 
Assume that $g\in\Cal G_K$ is generated by $f$ which satisfies
{\rm (3.12).} Let $K=[r\sqrt{\alpha}]$, $r>\sqrt{2C_0e^2}${\rm .}
Then{\rm ,}
$$
\eplana[g]=(1+O(r^{-K}))\frac 1{Z_K^\alpha}\sum_{h\in \Bbb
N^K}
\Delta_K(h)^2\prod_{i=1}^K\frac{\alpha^{h_i}}{h_i!^2}
\prod_{i=1}^Kf(h_i),\tag 3.16
$$
where
$$
Z_K^\alpha=\sum_{h\in \Bbb
N^K}
\Delta_K(h)^2\prod_{i=1}^K\frac{\alpha^{h_i}}{h_i!^2}.
$$
\endproclaim

\demo{Proof}
Write
$$
\eplana[g]=\eplana[g]-F_M^\alpha[g]+\frac{F_M^\alpha[1]}{\eplana[1]}
\frac{F_M^\alpha[g]}{F_K^\alpha[g]}
\frac{F_K^\alpha[1]}{F_M^\alpha[1]}
\frac{F_K^\alpha[g]}{F_K^\alpha[1]}.\tag 3.17
$$
By Lemma 3.2
$$
\lim_{M\to\infty}F_M^\alpha[g]=\eplana[g].\tag 3.18
$$
By Lemma 3.4 and the choice of $K$ and $r$
$$
\frac{F_K^\alpha[g]}{F_M^\alpha[g]}=1+O(r^{-K}),\tag 3.19
$$
for any $M>K$, and similarly with $g$ replaced by 1. Using (3.18) and
(3.19) in (3.17) and letting $M\to\infty$ we obtain
$$
\eplana[g]=(1+O(r^{-K}))\frac{F_K^\alpha[g]}{F_K^\alpha[1]},
$$
which is exactly (3.16). The proposition is proved.
\enddemo

Thus we have an approximate Coulomb gas picture of the (shifted) rows of
$\lambda$ under the Plancherel measure analogous to Dyson's Coulomb gas 
picture for the eigenvalues of a random matrix.

\demo{{R}emark {\rm 3.6}}  The
confining potential for the discrete Coulomb gas in (3.16) is
$$
V_K^\alpha[h_i]=-\frac 1K\log(\alpha^{h_i}/(h_i!)^2)
$$
with limit
$$
\lim_{\alpha\to\infty}V_K^\alpha[Kx]=2[x\log x+(\log r-1)x]=V(x).
$$
We can now use general techniques for Coulomb gases, see e.g. [Jo1],
[Jo3], to deduce asymptotic distribution properties. The potential $V$
has the (constrained) equilibrium measure; compare with Section 2 in [Jo3],
$u(t)dt$, where
$$
u(t)=\cases 1, &\text{if $0\le t\le 1-2/r$}\\
\frac 12-\frac 1{\pi}\arcsin(\frac r2(t-1)),&\text{if $1-2/r\le t\le 
1+2/r$}\\ 0,&\text{if $t\ge 1+2/r$.}\endcases
$$
Pick $f(t)=\exp(\phi(t/[r\sqrt{\alpha}]))$ with $\phi:\Bbb R\to\Bbb R$
continuous, bounded together with its derivative 
and $\phi(t)=0$ if $t\le 0$. Then,
$$
g(\lambda)=\prod_{i=1}^\infty\exp(\phi(\frac{\lambda_i+[r\sqrt{\alpha}]-i} 
{[r\sqrt{\alpha}]})).
$$
If we pick $r$ sufficiently large (depending on $\phi$) we can use
(3.19) and (3.20) to show that
$$
\lim_{\alpha\to\infty}\frac
1{[r\sqrt{\alpha}]}\log\eplana[g(\lambda)]=
\int_0^{1+2/r}\phi(t)u(t)dt.\tag 3.20
$$
From the limit (3.20) 
it is possible to deduce Vershik and Kerov's $\Omega$-law
for the asymptotic shape of the Young diagram, [VK], see also [AD],
where an outline of the argument using the hook-integral is
given. (The $r$-dependence in the formulas above goes away after
appropriate rescaling.) From what has been said above we see that the
$\Omega$-law is directly related to an equilibrium measure for a
discrete Coulomb gas. Using the general results in [Jo3] we can 
use (3.16) to show upper- and lower-tail large deviation formulas for
$\lambda_1$ ($=$ the length of the longest increasing subsequence in a
random permutation) under the Poissonized Plancherel measure. 
These formulas have been proved in [Se1], [Jo2] and [DZ] by other
methods
and we 
will not give the details of the new proof.
\enddemo
 
\demo{{\rm 3.3.} Proof of Theorem {\rm 1.3}}
We will now use Theorem 1.2 to prove Theorem 1.3, but before we can do
this we need certain asymptotic results for Charlier polynomials and
Bessel functions. Let
$$
w_a(x)=e^{-a}\frac{a^x}{x!},\quad x\in\Bbb N,\quad a>0.
$$
The normalized Charlier polynomials, $c_n(x;a)$, $n\ge 0$ are
orthonormal on $\Bbb N$ with respect to this weight.
The relevant value of the parameter $a$ for us will be $a=\alpha/M$,
and we define the {\it Charlier kernel}
$$\align
\KCh (x,y)=&\ \sqrt{\alpha}\frac{c_M(x;
\frac{\alpha}{M})c_{M-1}(y;\frac{\alpha}{M})
-c_{M-1}(x;\frac{\alpha}{M})c_M(y;\frac{\alpha}{M})}{x-y}\tag 3.21\\
&\times w_{\alpha/M}(x)^{1/2}w_{\alpha/M}(y)^{1/2},
\endalign$$
for $x\neq y$ and
$$ \multline
\KCh(x,x)\\
=\sqrt{\alpha}w_{\alpha/M}(x)\left[c_M'\left(x;
\frac{\alpha}{M}\right)c_{M-1}\left(x;\frac{\alpha}{M}\right)-c_{M-1}\left(x;\frac{\alpha}{M}\right)c_M'\left(x;
\frac{\alpha}{M}\right)\right].\endmultline\tag 3.22
$$
The polynomials $c_{n}(x;\alpha/M)$, $n\ge 0$, have the
generating function
$$
\sum_{n=0}^\infty\biggl(\frac{\alpha}M\biggr)^{n/2}\frac 1{\sqrt{n!}}
c_{n}\left(x;\frac{\alpha}{M}\right)w^n=e^{-\alpha w/M}(1+w)^x.
$$
It follows from this formula that we have the following integral
representations. If $0<r\le\sqrt{\alpha}/M$, then
$$
c_{n}(x;\frac{\alpha}{M})=\sqrt{\frac{n!}{M^n}}\frac
1{2\pi}\int_{-\pi}^\pi
e^{-\sqrt{\alpha}re^{i\theta}}\left(1+\frac{Mre^{i\theta}}
{\sqrt{\alpha}}\right)^x\frac 1{(re^{i\theta})^n}d\theta\tag 3.23
$$
and if $\sqrt{\alpha}/M<r$, then
$$\align
c_{n}(x;\frac{\alpha}{M})=&\ \sqrt{\frac{n!}{M^n}}\frac
1{2\pi}\int_{-\pi}^\pi
e^{-\sqrt{\alpha}re^{i\theta}}\left(1+\frac{Mre^{i\theta}}
{\sqrt{\alpha}}\right)^x\frac 1{(re^{i\theta})^n}d\theta\tag 3.24\\
&-(-1)^n\frac{\sin\pi x}{\pi}\int_{\sqrt{\alpha}/M}^r
e^{\sqrt{\alpha}s} \left(\frac{Ms}{\sqrt{\alpha}}-1\right)^xs^{-n}\frac
{ds}s,
\endalign
$$
for any $x\in\Bbb R$, where the powers are defined using the principal
branch of the logarithm.  

We want to write the Charlier kernel in a
form that will be convenient for later asymptotic analysis. Define,
for a given $r>0$, $x\in\Bbb Z$,
$$
A_M^\alpha(x)=\sqrt{\alpha}\frac {M!}{M^M}w_{\alpha/M}(x)\left(1+\frac{M}{
\sqrt{\alpha}}\right)^{2x}e^{-2\sqrt{\alpha}},
$$
$$
D_M^{\alpha,r}(x,g)=\frac1{2\pi}\int_{-\pi}^\pi g(re^{i\theta})
e^{\sqrt{\alpha}(1-re^{i\theta})}\left(\frac{\sqrt{\alpha}+Mre^{i\theta}}
{\sqrt{\alpha}+M}\right)^x\frac 1{(re^{i\theta})^M}d\theta,
$$
$$
F_M^{\alpha,r}(x,g)=(-1)^{x+M+1} \int_{\sqrt{\alpha}/M}^rg(s)
e^{\sqrt{\alpha}(1+s)}\biggl|\frac{\sqrt{\alpha}-Ms}{\sqrt{\alpha}+M}\biggr|^x 
s^{-M}\frac {ds}s,
$$
if $r>\sqrt{\alpha}/M$, and if
  $r\le\sqrt{\alpha}/M$, then $F_M^{\alpha,r}(x,g)=0$.
Then, some computation shows that when $x$ is an {\it integer} (the
case we are interested in),
$$
\KCh(x,y)=\sqrt{\!A_M^\alpha(x)A_M^\alpha(y)}
\frac{D_M^{\alpha,r}(x,g_1)D_M^{\alpha,r}(y,g_2)-D_M^{\alpha,r}(x,g_2)
D_M^{\alpha,r}(y,g_1)}{x-y},\tag 3.25
$$
when $x\neq y$, and
$$\align
\KCh(x,x)=&\ A_M^\alpha(x)\bigl[D_M^{\alpha,r}\!(x,g_2)D_M^{\alpha,r}(x-1,g_3)
-D_M^{\alpha,r}(x,g_1)D_M^{\alpha,r}(x-1,g_4)\bigr]\tag 3.26\\
&+A_M^\alpha(x)\bigl[F_M^{\alpha,r}(x,g_1)
D_M^{\alpha,r}(x,g_2)-F_M^{\alpha,r}(x,g_2)
D_M^{\alpha,r}(x,g_1)\bigr],
\endalign
$$
where $g_1(z)\equiv 1$, $g_2(z)=z-1$, 
$$
g_3(z)=
\biggl(\frac{\sqrt{\alpha}+Mz}{\sqrt{\alpha}+M}\biggr)\log
\biggl(\frac{\sqrt{\alpha}+Mz}{\sqrt{\alpha}+M}\biggr) ,
$$
and $g_4(z)=g_2(z)g_3(z)$. Note that all the $g_i$'s are bounded on
$|z|=r$.

The discrete Bessel kernel is defined by (1.10)
for $x\neq y$ and
$$
B^\alpha(x,x)=\sqrt{\alpha}[L_x(2\sqrt{\alpha})J_{x+1}(2\sqrt{\alpha})-
J_x(2\sqrt{\alpha})L_{x+1}(2\sqrt{\alpha})]\tag 3.27
$$
for $x=y$, where $L_x(t)=\frac {d}{dx}J_x(t)$.
The Bessel function has the integral representation
$$
J_x(2\sqrt{\alpha})=\frac 1{2\pi}\int_{-\pi}^\pi
e^{\sqrt{\alpha}(\frac
1re^{-i\theta}-re^{i\theta})+ix\theta}r^xd\theta-
\frac{\sin\pi x}{\pi}\int_0^re^{\sqrt{\alpha}(-1/s+s)}s^x\frac
{ds}s,\tag 3.28
$$
for $x\in\Bbb R$, $r>0$. Differentiation shows that for {\it integer}
$x$,
$$\align
L_x(2\sqrt{\alpha})=&\ \frac 1{2\pi}\int_{-\pi}^\pi
\log(re^{i\theta})e^{\sqrt{\alpha}(\frac
1re^{-i\theta}-re^{i\theta})+ix\theta}r^xd\theta\tag 3.29\\&-
(-1)^x\int_0^re^{\sqrt{\alpha}(-1/s+s)}s^x\frac
{ds}s.
\endalign$$
The next lemma shows that the discrete Bessel kernel is the
$M\to\infty$ limit of the Charlier kernel and establishes some
technical estimates. (We will only consider the case when $x,y$ are
integers but this restriction can be removed.)

\nonumproclaim{Lemma 3.7} For any $x,y\in\Bbb Z${\rm ,}

 \noindent {{\rm (i)}} 
$$\lim_{M\to\infty}\KCh(M+x,M+y)=B^\alpha(x,y).\tag 3.30$$

 \noindent {{\rm (ii)}} 
$$B^\alpha(x,y)=\sum_{k=1}^\infty J_{x+k}(2\sqrt{\alpha})
J_{y+k}(2\sqrt{\alpha}).\tag 3.31$$
  Furthermore{\rm ,} there is a constant $C=C(\alpha,L)${\rm ,} such that

 \noindent  {\rm (iii)}
$$
\sum_{x=-L}^\infty \KCh(M+x,M+x)\le C\tag 3.32
$$
 if $M$ is large enough{\rm ,} and

 \noindent {\rm (iv)} 
$$
\sum_{x=-L}^\infty B^\alpha(x,x)\le C.\tag 3.33
$$
 {\rm (}\/In {\rm (3.33)} we can take $C(\alpha,L)=\alpha/\sqrt{2}+L${\rm .)}
\endproclaim

\demo{Proof} We have to show that (3.25) and (3.26) converge to (1.10)
and (3.27) respectively. Using Stirling's formula we see that
$A_M^\alpha(M+x)\to\sqrt{\alpha}$ as $M\to\infty$. The result then
follows from the integral formulas above, the fact  that
$$
\lim_{M\to\infty}e^{\sqrt{\alpha}(1-z)}
\biggl(\frac{\sqrt{\alpha}+Mz}{\sqrt{\alpha}+M}\biggr)^{x+M}\frac
1{z^M}= e^{\sqrt{\alpha}(1/z-z)} z^y,
$$
and $g_3(z)\to z\log z$ as $M\to\infty$. This establishes (3.30). The
identity (3.31) follows from the recursion relation
$J_{x+1}(t)=2xJ_x(t)/t-J_{x-1}(t)$,
which implies
$$
B^\alpha(x,y)=J_{x+1}(t)J_{y+1}(t)+B^\alpha(x+1,y+1),
$$
and (3.31) follows by using the decay properties of the Bessel
function; see Lemma 3.9 below.

The estimate (3.32) is proved using the formula (3.26). Stirling's
formula can be used to show that $A_M^\alpha(x+M)\le 2\sqrt{\alpha}$
for all $x\ge 0$. We have
$$\align
\biggl|\frac{\sqrt{\alpha}+Mz}{\sqrt{\alpha}+M}\biggr|^M\frac 1{|z|^M}
\biggl|\frac{\sqrt{\alpha}+Mz}{\sqrt{\alpha}+M}\biggr|^y&\le
\left(1+\frac{\sqrt{\alpha}}{M|z|}\right)^M\left(|z|+\frac{\sqrt{\alpha}}{M}\right)^y\\
&\le \exp((1-\delta)^{-1}\sqrt{\alpha})(1-\delta/2)^y
\endalign
$$
if $|z|=r=1-\delta$ and $M\ge 2\sqrt{\alpha}/\delta$. This estimate
can be used in the integral formulas for $D_M^{\alpha,r}$ and 
$F_M^{\alpha,r}$ and we obtain
$$
|D_M^{\alpha,1/2}(M+x;g_i)|,|F_M^{\alpha,1/2}(M+x;g_i)|\le
Ce^{4\sqrt{\alpha}}(\frac 34)^x.
$$
Thus,
$$
\sum_{x=-L}^\infty \KCh(M+x,M+x)\le C\sqrt{\alpha}e^{4\sqrt{\alpha}}
\sum_{x=-L}^\infty
(\frac 34)^x.
$$

The estimate (iv) can be proved in a similar way but we can also
proceed as follows. Using the generating function for the Bessel
functions $J_n(t)$, $n\in\Bbb Z$, one can show that
$\sum_{n=1}^\infty n^2J_n(t)^2=t^2/4$,
see [Wa, 2.72(3)], and $\sum_{n=1}^\infty J_n(t)^2=\frac
12(1-J_0(t)^2)\le 1/2$, so by (ii) and the fact that $B(x,x)\le 1$,
$$
\sum_{x=-L}^\infty B^\alpha(x,x)\le L+\sum_{n=1}^\infty
nJ_n(2\sqrt{\alpha})^2\le \alpha/\sqrt{2}+L,
$$
where we have used the Cauchy-Schwarz  inequality. The lemma
is proved.
\enddemo

We are now ready for the
\demo{Proof of Theorem 1.3} We have
$$
\echara[g]=\prod_{j=1}^{M-1}\frac 1{j!}\sum_{\lambda\in\Omega_M}
\prod_{i=1}^Mf(\lambda_i+L-i)V_M(\lambda)^2W_M(\lambda)
\prod_{i=1}^M\biggl(\frac{\alpha}M\biggr)^{\lambda_i}e^{-\alpha/M}.
$$
If we make the change of variables $h_i=\lambda_i+M-i$, this can be
written
$$
\echara[g]=\frac 1{Z_M^\alpha}\sum_{h\in\Bbb N^M}
\prod_{i=1}^M(1+\phi(h_i-M+L))\Delta_M(h)^2\prod_{i=1}^Mw_{\alpha/M}(h_i).
$$
Now, using a standard computation from random matrix theory, see [Me],
[TW2], we can write this as
$$
\echara[g]
=\sum_{k=0}^{M-1}\sum_{h\in\Bbb N^k}
\prod_{i=1}^k\phi(h_i)\det(\KCh(h_i+M-L,h_j+M-L))_{i,j=1}^k\tag 3.34
$$
since $\phi(t)=0$ if $t<0$. The Charlier kernel is positive
definite, so we have the estimate
$$
\bigr|\det(\KCh(x_i,x_j))_{i,j=1}^k\bigr|\le\prod_{j=1}^k\KCh(x_j,x_j).
$$
Thus, by Lemma 3.7(iii),
$$
\multline
 \left|\sum_{h\in\Bbb N^k}
\prod_{i=1}^k\phi(h_i)\det(\KCh(h_i+M-L,h_j+M-L))_{i,j=1}^k\right|\\
 \le||\phi||_\infty^k\biggl(\sum_{x=-L}^\infty\KCh(M+x,M+x)\biggr)^k\le
(C||\phi||_\infty)^k.
\endmultline$$
The analogous estimate for the Bessel kernel follows from Lemma
3.7(iv). These estimates and Lemma 3.7(i) allow us to take the
$M\to \infty$
limit in (3.34). By Theorem 1.2 this gives (1.18). The theorem is
proved.
\enddemo

Note that we could just as well use Theorem 1.1 and the Meixner
ensemble to prove Theorem 1.3. The proof would be the same and we just
have to prove (3.30) and (3.32) for the Meixner ensemble instead.

\demo{{\rm 3.4.} Asymptotics of the Plancherel measure}
Theorem 1.3 can be used to analyze the asymptotic properties of the
Plancherel measure in different regions. One can distinguish three
cases corresponding to three different scaling limits of the Bessel
kernel. First we have the {\it edge scaling limit},
$$
\lim_{\alpha\to\infty} \alpha^{1/6}B^\alpha(2\sqrt{\alpha}+\xi
\alpha^{1/6},2\sqrt{\alpha}+\eta
\alpha^{1/6})=A(\xi,\eta),\tag 3.35
$$
where $A$ is the {\it Airy kernel} defined in (1.4). This is the case that
is considered in Theorem 1.4. Secondly we
have the {\it bulk scaling limit},
$$
\lim_{\alpha\to\infty}
B^\alpha(r\sqrt{\alpha},r\sqrt{\alpha}+u)=\frac{\sin(uR)}{u\pi},\tag
3.36
$$
$u\in\Bbb Z$, $-2<r<2$, where $R=\arccos(r/2)$; the right-hand side is
the {\it discrete sine kernel}. We will not discuss the local
behavior in the bulk of the Young diagram; see [BOO]. Thirdly we have
an {\it intermediate region},
$$\multline
\lim_{\alpha\to\infty}\pi\alpha^{1/4-\delta/2}B^\alpha
(2\sqrt{\alpha}-\alpha^\delta+\pi\xi\alpha^{1/4-\delta/2}
,2\sqrt{\alpha}-\alpha^\delta+\pi\eta\alpha^{1/4-\delta/2})\\
 =
\frac{\sin\pi(\xi-\eta)}{\pi(\xi-\eta)},\endmultline \tag 3.37
$$
if $1/6<\delta<1/2$, the ordinary {\it sine kernel}. Thus in this region the
local behavior is the same as that in the bulk in a random Hermitian
matrix. The limits (3.35) to (3.37) can be proved using the
saddle-point method on the integral formula for the Bessel
function. From the point of view of the Coulomb gas picture of the
Young diagram, the cases one and three are similar to the random
matrix case since at the edge a discrete Coulomb gas approximates a
continuous Coulomb gas. Case two is different however, since in the
bulk the discrete nature is manifest; the charges sit close to each
other. 

Before turning to the proof of Theorem 1.4 we have to say something
about de-Poissonization, the joint distribution of the first $k$ rows
($k$ largest eigenvalues) and the asymptotics of the Bessel kernel.

We have the following generalization of a lemma in [Jo2].

\nonumproclaim{Lemma 3.8} Let $\mu_N=N+4\sqrt{N\log N}$ and
$\nu_N=N-4\sqrt{N\log N}$. Then there is a constant $C$ such that{\rm ,} for
$0\le x_i\le N${\rm ,}
$$\multline\Bbb P_{\text{Plan}}^{\mu_N}[\lambda_1\le x_1,\dots,\lambda_k\le
x_k]-\frac C{N^2}\le
\pplanchN[\lambda_1\le x_1,\dots,\lambda_k\le
x_k]\\
  \le
\Bbb P_{\text{Plan}}^{\nu_N}[\lambda_1\le x_1,\dots,\lambda_k\le
x_k]+\frac C{N^2}.
\endmultline\tag 3.38$$
\endproclaim

\demo{Proof} This is proved as Lemmas 2.4 and 2.5 in [Jo2]. Denote a
permutation in $S_N$ by $\pi^{(N)}$ and let $S_{N+1}(j)$ denote the
set of all $\pi^{(N+1)}$ such that $\pi^{(N+1)}(N+1)=j$. Each
$\pi^{(N+1)}$ in $S_{N+1}(j)$ is mapped to a permutation
$F_j(\pi^{(N+1)})$ in $S_N$ by replacing each $\pi^{(N+1)}(i)>j$ by 
$\pi^{(N+1)}(i)-1$. The map $F_j$ is a bijection from $S_{N+1}(j)$ to
$S_N$. Apply the Robinson-Schensted correspondence to 
$F_j(\pi^{(N+1)})$ to obtain the $P$-tableau. Replace the entries $i$ by
$i+1$ for $i=j,\dots, N$ and then insert $j$. This insertion can only
increase the length of any row and we obtain the $P$-tableau for
$\pi^{(N+1)}$. Thus,
$$
\lambda_i(F_j(\pi^{(N+1)}))\le\lambda_i(\pi^{(N+1)}),
$$
for all rows. If we define $g(\pi^{(N)})$ to be 1 if
$\lambda_i(\pi^{N})\le x_i$ for $i=1,\dots,k$ and 0 otherwise, we see
that
$$
g(F_j(\pi^{(N+1)}))\ge g(\pi^{(N+1)}),
$$
and we can proceed exactly as in [Jo2] using the fact that the
Plancherel measure on $\Cal P^{(N)}$ is the push-forward of the
uniform distribution on $S_N$.
\enddemo

For $x\in\Bbb R^M$, $n\in\Bbb N^k$ and a sequence $\Cal I= (I_1,\dots,
I_k)$ of intervals in $\Bbb R$ we let $\chi(\Cal I,n,x)$ denote the
characteristic function for the set of all $x\in\Bbb R^M$ such that
exactly $n_j$ of the $x_i$'s belong to $I_j$, $j=1,\dots,k$. A
computation shows that for a single interval
$$
\chi(I_j,n_j,x)=\frac 1{n_j!}\frac{\partial^{n_j}}{\partial z_j^{n_j}}
\prod_{i=1}^M(1+z_j\chi_{I_j}(x_i))\biggr|_{z_j=-1}
$$
and hence
$$
\chi(\Cal I,n,x)=\frac 1{n_1!\dots n_k!}
\frac{\partial^{n_1+\dots+n_k}}{\partial z_j^{n_1}\dots\partial
z_j^{n_k}}\prod_{i=1}^M\prod_{j=1}^k(1+z_j\chi_{I_j}(x_i)) 
\biggr|_{z_1=\dots=z_k=-1}.\tag 3.39
$$
Note that if the intervals are pairwise disjoint, then 
$\prod_{j=1}^k(1+z_j\chi_{I_j}(x_i))=1+\sum_{j=1}^k
z_j\chi_{I_j}(x_i)$;
compare with [TW2]. Let $\Bbb P$ be a probability measure on $\Bbb R^M$ and
let $a_1\ge\dots\ge a_k$. Set $I_{j+1}=(a_{j+1},a_j]$, $j=1,\dots,
k-1$ and $I_1=(a_1,\infty)$. Let
$$
\Bbb L_k=\{n\in\Bbb N^k\,;\,\sum_{j=1}^rn_j\le r-1, \,r=1,\dots,k\}.
$$
Define $x^{(j)}$ to be the $j^{\rm th}$ largest of the $x_i$'s. Then,
$$
\Bbb P[x^{(1)}\le a_1,\dots,x^{(k)}\le a_k]=\sum_{n\in \Bbb L_k}
\Bbb E[\chi(\Cal I,n,x)].\tag 3.40
$$
Hence, the problem of investigating the distribution function in
(3.40) reduces to investigating expectations of the right-hand side of
(3.39). 

In the proof of Theorem 1.4 we will need some asymptotic results for
Bessel functions. 

\nonumproclaim{Lemma 3.9} Let $M_0>0$ be given. Then there exists a
constant $C=C(M_0)$ such that if we write
$x=2\sqrt{\alpha}+\xi\alpha^{1/6}${\rm ,} then
$$
|J_x(2\sqrt{\alpha})|\le C\alpha^{-1/6}\exp\left[-\frac 14\min\left(\frac
14\alpha^{1/6}, |\xi|^{1/2}\right)|\xi|\right]\tag 3.41
$$
for $\xi\in [-M_0,\infty)$. Furthermore
$$
\lim_{\alpha\to\infty}\alpha^{1/6}J_x(2\sqrt{\alpha})=\Ai(\xi),\tag
3.42
$$
uniformly for $\xi\in[-M_0,M_0]$.
\endproclaim

This can be deduced from classical asymptotic results, [Wa]
and it is also rather straightforward to proceed as in Section 5 of
[Jo3] using the integral formula (3.28).

We are now ready for the
\demo{Proof of Theorem 1.4} We will prove (1.20).  The proof of (1.19)
is analogous using the Hermite kernel instead. From Lemma 3.8, the
fact that a distribution function is increasing in its arguments, that
the distribution function $F(t_1,\dots,t_k)$ is continuous and 
$\sqrt{\mu_N}-\sqrt{N}\approx 2\sqrt{\log N}$, 
$\sqrt{\nu_N}-\sqrt{N}\approx 2\sqrt{\log N}$, we see that it suffices
to prove that
$$
\lim_{\alpha\to\infty}\pplancha[\lambda_1-1\le
2\sqrt{\alpha}+t_1\alpha^{1/6}, \dots
\lambda_k-k\le
2\sqrt{\alpha}+t_k\alpha^{1/6}]=F(t_1,\dots,t_k),\tag 3.43
$$
for any fixed $(t_1,\ldots,t_k)\in\Bbb R^k$, $t_1\ge\dots\ge t_k$. Set $$I_{j+1}=(
2\sqrt{\alpha}+t_{j+1}\alpha^{1/6},2\sqrt{\alpha}+t_j\alpha^{1/6}],\qquad j=1,\dots,k-1$$
 and $I_1=(2\sqrt{\alpha}+t_1\alpha^{1/6},\infty)$. By
(3.39) and (3.40) it is enough to consider the expectations
$$
\eplana\left[\prod_{i=1}^\infty\prod_{j=1}^k(1+z_j\chi_{I_j}(\lambda_j-j))\right].
\tag 3.44
$$
If we write $\phi_\alpha(s)=\prod_{j=1}^k(1+z_j\chi_{I_j}(s))-1$ it
follows from Theorem 1.3, with $L=0$, that the expectation (3.44) can
be written as
$$
F_\alpha(z,t)=\sum_{k=0}^\infty\frac 1{k!}\sum_{h\in\Bbb N^k}
\prod_{j=1}^k\phi_\alpha(h_j)\det[B^\alpha(h_i,h_j)]_{i,j=1}^k.\tag
3.45
$$
Note that $F_\alpha(z,t)$ is an entire function of $z$. Set
$J_{j+1}=(t_{j+1},t_j]$, $j=1,\dots,\break k-1$, $J_1=(t_1,\infty)$ and write
$\psi(s)=\prod_{j=1}^k(1+z_j\chi_{J_j}(s))-1$. Define
$$
F(z,t)=\sum_{k=0}^\infty\frac 1{k!}\int_{\Bbb R^k}\prod_{j=1}^k
\psi(\xi_j)\det[A(\xi_i,\xi_j)]_{i,j=1}^k.\tag 3.46
$$
We want to show that
$$
\lim_{\alpha\to\infty}F_\alpha(z,t)=F(z,t),\tag 3.47
$$ 
uniformly for $z$ in a
compact subset of $\Bbb C^k$. Then also derivatives of $F_\alpha(z,t)$
converge to the corresponding derivatives of $F(z,t)$. The limit
(3.43) then follows with
$$
F(t_1,\dots,t_k)=\sum_{n\in\Bbb L_k}\frac 1{n_1!\dots n_k!}
\frac{\partial^{n_1+\dots +n_k}}{\partial z_1^{n_1}\dots
\partial z_k^{n_k}}F(z,t)\biggr|_{z_1=\dots z_k=-1}.\tag 3.48
$$

So it remains to prove (3.47). Note that $\phi_\alpha(s)=0$ if
$s<2\sqrt{\alpha} +t_k\alpha^{1/6}$ and that $\phi_{\alpha}(s)=
\psi(\alpha^{-1/6}(s-2\sqrt{\alpha}))$. Given $r\in\Bbb R$ we set
$\Bbb A(r)=\{r,r+1,r+2,\dots\}$.  Then,
$$
F_\alpha(z,t)=\sum_{l=0}^\infty\frac 1{l!}\sum_{h\in\Bbb
A(t_k\alpha^{1/6})^l} \prod_{j=1}^l\psi\left(\frac{h_j}{\alpha^{1/6}}\right)
\det[\tilde B^\alpha(\xi,\eta)]\frac 1{(\alpha^{1/6})^l},\tag 3.49
$$
where $\tilde B^\alpha(\xi,\eta)
=\alpha^{1/6}B^\alpha(2\sqrt{\alpha}
+\xi\alpha^{1/6},2\sqrt{\alpha} +\eta\alpha^{1/6})$.
We can now prove that (3.47) holds pointwise in $z$ by the same
argument as was used in the proof of the analogous statement in
Section 3 of [Jo3]. That proof depends on the following properties of
the kernel; compare with Lemma 3.1 in [Jo3] and Lemma 4.1 below.
\item{(i)} For any $M_0>0$ there is a constant $C=C(M_0)$ such that
for all $\xi\ge -M_0$
$$
\sum_{m=1}^\infty B^\alpha(2\sqrt{\alpha}
+\xi\alpha^{1/6}+m,2\sqrt{\alpha}
+\xi\alpha^{1/6}+m)\le C.
$$
\item{(ii)} For any $\varepsilon>0$, there is an $L>0$ such that
$$
\sum_{m=1}^\infty B^\alpha(2\sqrt{\alpha}
+L\alpha^{1/6}+m,2\sqrt{\alpha}
+L\alpha^{1/6}+m)\le\varepsilon,
$$
for all sufficiently large $\alpha$.
\item{(iii)} For any $M_0>0$ and any $\varepsilon>0$
$$
\left|\tilde B^\alpha\left(\frac n{\alpha^{1/6}},\frac m{\alpha^{1/6}}\right)-
A\left(\frac n{\alpha^{1/6}},\frac m{\alpha^{1/6}}\right)\right|\le\varepsilon
$$
for all integers $m,n\in[-M_0\alpha^{1/6},M_0\alpha^{1/6}]$ provided
$\alpha$ is sufficiently large.

The estimate (i) is used to estimate the tail in the $k$-summation in
(3.49), (ii) is used to limit the $h$-summation and (iii) is used to
prove that the Riemann sums converge to integrals.

If $z$ belongs to a compact set $K$ there is a constant $C$,
independent of $z$, such that $||\psi||_\infty\le C$. Together with
(i) this shows that the family $\{F_\alpha(z,t)\}$ is uniformly
bounded for $\alpha>0$, $z\in K$ and hence (3.47) holds uniformly
by a normal family argument.

The properties (i) to (iii) above are straightforward to prove using
the representation (3.31) and Lemma 3.9. To prove (i) and (ii) we use
$$
\sum_{m=1}^\infty B^\alpha(x+m,x+m)=\sum_{n=1}^\infty
nJ_{x+n+1}^2(2\sqrt{\alpha}),
$$
which can be estimated using (3.41) (we get a Riemann sum). Similarly,
\linebreak 
$\tilde B^\alpha(\frac n{\alpha^{1/6}},\frac m{\alpha^{1/6}})$ can be
written as a Riemann sum, using (3.31), which is controlled using
(3.41) and (3.42). This Riemann sum can be compared with the
corresponding Riemann sum for the following representation of the Airy
kernel, [TW1],
$$
A(\xi,\eta)=\int_0^\infty \Ai(\xi+t)\Ai(\eta+t)dt
$$
and in this way we obtain (iii).

\section{Random words and the Charlier ensemble}

In this section we will prove our results on random words.

\demo{Proof of Proposition {\rm 1.5}} Let $L(M,N,\lambda)$ denote the number
of pairs\break $(P,Q)$ of tableaux of shape $\lambda\in\Omega_M^{(N)}$ with
$P$ semistandard with elements in $\{1,\dots,M\}$ and $Q$ standard
with elements in $\{1,\dots,N\}$. Then
$$
\pwMN[S^{-1}(\lambda)]=\frac 1{M^N}L(M,N,\lambda).\tag 4.1
$$
The number of possible $P$'s is, by [Fu],
$$
d_\lambda(M)=\prod_{1\le i< j\le M}\frac{\lambda_i-\lambda_j+j-i}{j-i}=
\left(\prod_{j=1}^{M-1}\frac 1{j!}\right)V_M(\lambda),\tag 4.2
$$
and the number of possible $Q$'s is $f^\lambda$ given by (1.13).
By (4.2), (4.3) and Lemma 3.1 we obtain
$$
L(M,N,\lambda)=N!\biggl(\prod_{j=1}^{M-1}\frac
1{j!}\biggr)V_M(\lambda)^2W_M(\lambda).
\tag 4.3
$$
Inserting the formula (4.3) into (4.1) yields the desired result
(1.21). The formulas (1.22) and (1.23) are immediate consequences. The
proposition is proved.
\enddemo
Next, we give the
\demo{Proof of Theorem 1.6} We will prove (1.24); the proof of (1.25)
is analogous. Both are straightforward asymptotic computations using
Stirling's formula and we will indicate the main steps. Set
$$
x_j=\frac{\lambda_j-N/M}{\sqrt{2N/M}},\quad j=1,\dots,M.
$$
Note that $\sum_{j=1}^M x_j=0$, since $\sum_{j=1}^M\lambda_j=N$. Then,
$$
(\lambda_j+M-j)!=\sqrt{\frac{2\pi N}{M}}\biggl(\frac
NM\biggr)^{N/M+M-j} e^{x_j^2-N/M+o(1)}
$$
as $N\to\infty$, and hence
$$
W_M(\lambda)\sim\biggl(\frac{2\pi N}{M}\biggr)^{-M/2}e^N
\biggl(\frac MN\biggr)^{N+M(M-1)/2}\prod_{j=1}^Me^{-x_j^2}.
$$
Furthermore,
$$
V_M(\lambda)^2=\biggl(\frac{2N}{M}\biggr)^{M(M-1)/2}\prod_{1\le
i<j\le M}\left(x_i-x_j+\frac{i-j}{\sqrt{2N/M}}\right),
$$
and consequently
$$\align
\pcharN[\{\lambda\}]&\sim\sqrt{\pi M}(2\pi)^{-M/2} 
2^{M^2/2}\prod_{j=1}^{M-1}\frac 1{j!}\Delta_M(x)^2\prod_{j=1}^Me^{-x_j^2}
\biggl(\frac{2N}{M}\biggr)^{-(M-1)/2}\tag 4.4\\
&=\sqrt{\pi M}M!\phi_{{\rm GUE},M}(x).
\endalign
$$
From this we see that the left-hand side of (1.24) is approximately a
Riemann sum for the right-hand side, which in the limit $N\to\infty$
converges to the right-hand side. The factor $M!$ in the last
expression in (4.4) comes from the fact that in (4.4) the variables
are ordered. This completes the proof.
\enddemo

For the proof of Theorem 1.7 we need asymptotic results for the
Charlier kernel analogous to those for the Bessel kernel in the proof
of Theorem 1.4.

\nonumproclaim{Lemma 4.1} Let $\nu=M+\alpha/M+2\sqrt{\alpha}$ and
$\sigma=(1+\sqrt{\alpha}/M)^{2/3}\alpha^{1/6}${\rm .}
\smallbreak
{{\rm (i)}} For any $M_0>0$ there is a constant $C=C(M_0)$ such that{\rm ,}
for all integers $n\ge -M_0\sigma${\rm ,}
$$
\sum^\infty_{m=1}\KCh([\nu]+n+m,[\nu]+n+m)\le C.\tag 4.5
$$
\smallbreak
{\rm (ii)} For any $\varepsilon>0$ there is an $L>0$ such that 
$$
\sum^\infty_{m=1}\KCh([\nu]+[\sigma L]+m,[\nu]+[\sigma L]+m)\le \varepsilon  \tag 4.6
$$
if $M,\alpha$ are sufficiently large{\rm .}
\medbreak
{\rm (iii)} For any $M_0>0$ and any $\varepsilon>0${\rm ,}
$$
\left|\sigma\KCh([\nu]+m,[\nu]+n)-A\left(\frac {m}{\sigma},\frac
{n}{\sigma}\right)\right|\le\varepsilon \tag 4.7
$$
for all integers $m,n\in[-M_0\sigma,M_0\sigma]$ provided $\alpha$ and
$M$ are sufficiently large{\rm .}
\endproclaim

\demo{Proof} 
The proof is based on the formulas (3.25) and (3.26) for the Charlier
kernel. The proof is completely analogous to the proof of the
corresponding result for the Meixner kernel in Lemma 3.2 in [Jo3, \S 5], so we will not give the details here. Asymptotic
formulas for Charlier polynomials with fixed $a=\alpha/M$ have been obtained in [Go].
\enddemo

\demo{Proof of Theorem {\rm 1.7}} By (1.16) and (1.23)
$$
\Bbb P_{\text{W},M}^\alpha[L(w)\le s]=
\prod_{j=1}^M\frac 1{j!}\sum\Sb h\in\Bbb N^M\\\max h_j\le s+M-1\endSb
\Delta_M(h)^2\prod_{j=1}^Mw_{\alpha/M}(h_i),
$$
where we have made the substitution $h_i=\lambda_i+M-i$. Using 
Lemma 4.1 this can be analyzed exactly as the analogous problem involving
the Meixner weight in Section 3 in [Jo3]. Lemma 3.1 in [Jo3] gives
$$
\Bbb P_{\text{W},M}^\alpha\left[L(w)\le\frac{\alpha}M
+2\sqrt{\alpha}+\biggl(1+\frac
{\sqrt{\alpha}}M\biggl)^{2/3}\alpha^{1/6}\xi\right] \to F(\xi),\tag 4.8
$$
as $\alpha,M\to\infty$ with $F(\xi)$ given by (1.5). This proves
(1.26). Next, we observe that for fixed $M$,  
$\pwMN[L(w)\le s]$ is a decreasing function of $N$, which can be proved
as the corresponding result for permutations in [Jo2]. Thus, with
$\mu_N$ and $\nu_N$ as in Lemma 3.8, we have
$$
\Bbb P_{\text{W},M}^{\mu_N}[L(w)\le s]-\frac C{N^2}\le
\pwMN[L(w)\le s]\le\Bbb P_{\text{W},M}^{\nu_N}[L(w)\le s]+\frac
C{N^2}. \tag 4.9
$$
Set $s(\alpha,M,\xi)=\frac{\alpha}M
+2\sqrt{\alpha}+\bigl(1+\frac
{\sqrt{\alpha}}M\bigl)^{2/3}\alpha^{1/6}\xi$. Then, $s(N,M,\xi)=
s(\mu_N,M,\break\xi+\delta)$ and $s(N,M,\xi)=
s(\nu_N,M,\xi+\delta')$, where $\delta,\delta'\to 0$ as $M,N\to\infty$
if $M^{-1}(\log N)^{1/6}$ converges to $0$ as $M,N\to\infty$. Thus, (1.27)
follows from (4.8) and (4.9) and the theorem is proved.
\enddemo

\section{Applications of the Krawtchouk ensemble}

 5.1. {\it Sepp{\rm \"{\it a}}l{\rm \"{\it a}}inen\/{\rm '}\/s first passage percolation model}.
The {\it Krawtchouk
ensemble} is defined by (1.7) with the weight  
$w(x)=\binom{K}{x}p^xq^{K-x}$,
$0\le x\le K$, i.e. we consider the probability measure
%\medbreak \centerline{${\displaystyle
$$\Bbb P_{\text{Kr},N,K,p}[h]=\frac 1{Z_{N,K,p}}\Delta_N(h)^2
\prod_{j=1}^N\binom{K}{h_j}p^{h_j}q^{K-h_j}
$$ %}$}
%\medbreak\noindent 
on  $\{0,\dots,K\}^N$, where
$
Z_{N,K,p}=N!\biggl(\prod_{j=0}^{N-1}
\frac{j!}{(K-j)!}\biggr)(K!)^N(pq)^{N(N-1)/2}.
$
The first problem where the Krawtchouk ensemble appears is in 
the simplified first-passage percolation model 
introduced by Sepp\"al\"ainen in [Se2]. 
Consider the lattice $\Bbb N^2$ and attach a \it passage time \rm
$\tau(e)$ to each nearest neighbour edge. If $e$ is vertical $\tau(e)=\tau_0
>0$, and if $e$ is horizontal then $\tau(e)$ is random with $P[\tau(e)=\lambda]
=p$ and $P[\tau(e)=\kappa]=q=1-p$, where $\kappa>\lambda\ge 0$, $0<p<1$. All 
passage times assigned to horizontal edges are independent random variables. 
Hence, all randomness sits in the horizontal edges. The \it minimal
passage time \rm from $(0,0)$ to $(k,l)$ along nearest neighbour paths is 
defined by
$$
T(k,l)=\min_{p}\sum_{e\in p}\tau(e)
\tag 5.1
$$
where the minimum is over all
non-decreasing nearest
neighbour paths $p$  from $(0,0)$ to $(k,l)$.
The \it time constant \rm is defined by
$
\mu(x,y)=\lim_{n\to\infty}\!\frac 1{n}T([nx],[ny]).
$
(The existence of the limit follows from subadditivity.) In [Se2] it is proved,
using a certain associated stochastic process, that
$$
\mu(x,y)=\cases \lambda x+\tau_0y, &\text{if $py>qx$}\\
\lambda x+\tau_0 y+(\kappa-\lambda)(\sqrt{qx}-\sqrt{py})^2, 
&\text{if $py\le qx$.}\endcases\tag 5.2
$$
We will show that the distribution of the random variable $T(k,l)$ relates 
to the distribution of the rightmost charge (``largest eigenvalue'') in a
Krawtchouk ensemble.

Write $M=k, N=l+1$ and consider an $M\times N$ matrix $W$ whose elements, 
$w(i,j)$, are independent Bernoulli random variables, $P[w(i,j)=0]=q$ and
$P[w(i,j)=1]=p=1-q$. Let $\Pi_{M,N}$ be the set of all sequences $\pi=
\{(k,j_k)\}_{k=1}^M$ such that $1\le j_1\le\dots\le j_M\le N$, i.e. 
up/right paths in $W$ with exactly one
element in each row. Introduce the
random variable
$$
L(W)=\max\left\{\sum_{(i,j)\in\pi} w(i,j)\,;\,\pi\in\Pi_{M,N}\right\}.\tag 5.3
$$
Write $\rho=1/q-1$, so that $q=(1+\rho)^{-1}$ and $p=\rho(1+\rho)^{-1}$. 
It is straightforward to show that
$$
T(k,l)=l\tau_0+k\kappa-(\kappa-\lambda)L(W).\tag 5.4
$$
\nonumproclaim{Proposition 5.1} Let $L(W)$ be defined by {\rm (5.3)} with $W$ an $M\times
N$ $0-1$\/{\rm -}\/matrix with independent Bernoulli elements $w(i,j)${\rm ,} 
the probability of $1$
being~$p${\rm .} Then{\rm ,}
$$
P[L(W)\le n]=\Bbb P_{\text{Kr},N,N+M-1,p}\left[\max_{1\le j\le N}h_j\le
n+N-1\right].\tag 5.5
$$
\endproclaim

\demo{Proof}
Interpreting the formula (7.30) in Theorem 7.1 in [BR1] in the appropriate 
case, we get
$$
P[L(W)\le n]=(1+\rho)^{-MN}\sum\Sb\lambda\in\Cal P\\ \ell(\lambda)\le n
\endSb
d_\lambda(M)d_{\lambda'}(N)\prod_{i=1}^{\ell(\lambda)}\rho_{\lambda_i},
\tag 5.6
$$
where $\lambda'$ is the partition conjugate to $\lambda$, $\lambda_k'$ is 
the length of the $k^{\rm th}$ column in $\lambda$, and $d_\lambda(M)$ is the 
number of semi-standard tableaux of shape $\lambda$ with elements in $\{
1,\dots,M\}$; if $\ell(\lambda)\le M$, $d_\lambda(M)$ is given by (4.2).
The proof of (5.6) is based on the RSK-correspondence between 0-1 matrices
and pairs of semistandard Young tableaux $(P,Q)$ where $P$ has shape
$\lambda$ and $Q$ has shape $\lambda'$, see [Fu], [St].
Set 
$$
\Omega_M(N)=\{\lambda\in\Omega_M\,;\,N\ge\lambda_1\ge\dots\lambda_M\ge 0\}.
$$
Since $d_\lambda(M)=0$ if $\ell(\lambda)>M$ and $d_{\lambda'}(N)=0$ if
$\lambda_1>N$, (5.6) can be written as
$$ \multline
P[L(W)\le n]\\ =(1+\rho)^{-MN}\left(\prod_{j=1}^{M-1}\frac 1{j!}\right)
\left(\prod_{j=1}^{N-1}\frac 1{j!}\right)\sum\Sb \lambda\in\Omega_M(N)\\
\ell(\lambda)\le n\endSb V_M(\lambda)V_N(\lambda')\prod_{i=1}^M\rho^{\lambda
_i}.\endmultline\tag 5.7
$$
Note that $\lambda\in\Omega_M(N)$ if and only if $\lambda'\in\Omega_N(M)$ 
and $\ell(\lambda)=\lambda_1'$.

\nonumproclaim{Lemma 5.2} If $\mu\in\Omega_N(M)${\rm ,} then
$$
V_M(\mu')=\left(\prod_{j=1}^{N+M-1}j!\right)V_N(\mu)W_N(\mu)\prod_{j=1}^N
\frac 1{(M+j-1-\mu_j)!}.\tag 5.8
$$
\endproclaim

\demo{Proof} One way to prove (5.8) is to use the fact that $V_M(\mu')W_M(\mu')
$ \linebreak $=V_M(\mu)W_M(\mu)$ by the hook formula for $f^\mu$; 
compare with (1.13) and Lemma 3.1.
We will give another 
proof. Set $s_i=\mu_i+N+1-i$, $1\le i\le N$ and $r_j=N+j-\mu_j'$, 
$1\le j\le M$. Then,
$$
\{s_1,\dots,s_N\}\cup\{r_1,\dots,r_M\}=\{1,\dots,N+M\}.\tag 5.9
$$
To see this, notice that since $1\le s_i,r_j\le N+M$ it suffices to show that
$s_i\neq r_j$ for all $i,j$. Looking at the $\mu$-diagram one sees that
$\mu_i+\mu_j'\le i+j-2$ or $\mu_i+\mu_j'\ge i+j$, which implies 
$s_i\neq r_j$.

Let $n_k=1$ if $k\in\{s_1,\dots,s_N\}$ and $n_k=0$ if $k\in\{r_1,\dots,r_M\}$,
$k=1,\dots,N+M$. Then, by (5.9),
$$
V_M(\mu')=\prod_{1\le k< l\le N+M}(l-k)^{(1-n_k)(1-n_l)}.\tag 5.10
$$
Now,
$$\align
\prod_{1\le k<l\le N+M}(l-k)^{n_kn_l}&=V_N(\mu),\\
\prod_{1\le k<l\le N+M}(l-k)^{n_k}&=\prod_{j=1}^N\prod_{l=s_j+1}^{N+M}
(l-s_j)=\prod_{j=1}^N(N+M-s_j)!
\\
\noalign{\noindent 
and\hfill}
\prod_{1\le k<l\le N+M}(l-k)^{n_l}&=\prod_{j=1}^N\prod_{k=1}^{s_j-1}
(s_j-k)=\prod_{j=1}^N(s_j-1)!.\endalign
$$
Inserting this into (5.10) gives the formula (5.8). The lemma is proved.
\enddemo
 
We can now finish the proof of the proposition.
If we write $\mu=\lambda'$, we see from (5.8) that (5.7) can be written as
$$\align
P[L(W)\le n]=&\ (1+\rho)^{-MN}\prod_{j=0}^{N-1}\frac{(j+M)!}{j!} \tag 5.11
\\&\times\sum\Sb \mu\in\Omega_N(M)\\ \mu_1\le n\endSb V_N(\mu)^2W_N(\mu)
\prod_{j=1}^N\frac{\rho^{\mu_j}}{(M+j-1-\mu_j)!}.
\endalign$$
As usual we introduce the new coordinates $h_j=\mu_j+N-j$. Then, using
$\rho=1/q-1$, we obtain
$$\align
P[L(W)\le n]=&
\ \frac 1{N!}\prod_{j=0}^{N-1}\frac{(j+M)!}{j!}\frac
{(pq)^{N(N-1)/2}}{((N+M-1)!)^N}\\ &\times
\sum\Sb h\in\Bbb N^N\\ \max(h_j)\le n+N-1
\endSb
\Delta_N(h)^2\prod_{j=1}^N\binom{N+M-1}{h_j}p^{h_j}q^{N+M-1-h_j},
\endalign$$
which completes the proof.\hfill\qed
\enddemo

Using Proposition 5.1 we can prove a limit theorem for the first
passage time $T(k,l)$. The result should be compared with Remark 1.8
and Conjecture 1.9 in [Jo3].

\nonumproclaim{Theorem 5.3} If $\mu(x,y)$ is given by {\rm (5.2),}
$$
\sigma(x,y)=\frac{(pq)^{1/6}}{(xy)^{1/6}}(\sqrt{px}+\sqrt{qy})^{2/3}
(\sqrt{qx}-\sqrt{py})^{2/3}
$$
and $py<qx${\rm ,} then
$$\align
&\lim_{n\to\infty}\Bbb
P\bigl[\frac{T([nx],[ny])-n\mu(x,y)}{\sigma(x,y)n^{1/3}}\le\xi\bigr] 
=1-F(-\xi),
\\ \noalign{\noindent 
where $F(t)$ is the Tracy\/{\rm -}\/Widom distribution {\rm (1.5).}}\endalign$$
\endproclaim

\demo{Proof} The proof uses (5.4) and Proposition 5.1 and 
is analogous to the proof of Theorem 1.7, the
difference being that we now need the analogue of Lemma 4.1 for the
Krawtchouk polynomials. This can be obtained from a steepest descent
analysis of the integral
formula for these polynomials in much the same way as in the
analysis of the Meixner polynomials in Section 5 of [Jo3]; see [Jo4] for
some more details. The time constant is related to the right endpoint
of the support of the equilibrium measure associated with the
Krawtchouk ensemble, and the constant $\sigma(x,y)$ comes out of the
steepest descent argument. We can also get large deviation results by
using the general results of Section 4 in [Jo3].
\enddemo

\demo{{\rm 5.2.} The Aztec diamond}
We turn now to 
the relation between the\break Krawtchouk ensemble and
domino tilings of the Aztec diamond introduced by Elkies, 
Kuperberg, Larsen and Propp in [EKLP]. The definitions are taken from that 
paper and the papers [JPS] and [CEP] where more details and pictures can be 
found. A \it domino \rm is a closed $1\times 2$ or $2\times 1$ rectangle 
in $\Bbb R^2$ with corners in $\Bbb Z^2$, and a \it tiling \rm of a region 
$R\subseteq\Bbb R^2$ by dominoes is a set of dominoes whose interiors are
disjoint and whose union is $R$. The \it Aztec diamond\rm, $A_n$, of order $n$
is the union of all lattice squares $[m,m+1]\times[l,l+1]$, $m,l\in\Bbb Z$,
that lie inside the region $\{(x,y)\,;\, |x|+|y|\le n+1\}$. It is proved in 
[EKLP] that the number of possible domino tilings of $A_n$ equals 
$2^{n(n+1)/2}$. Color the Aztec diamond in a checkerboard fashion so that
the leftmost square in each row in the top half is white. A horizontal domino
is \it north-going \rm if its leftmost square is white, otherwise it is 
\it south-going\rm. Similarly, a vertical domino is \it west-going \rm if its 
upper square is white, otherwise it is \it east-going\rm. Two dominoes are
\it adjacent \rm 
if they share an edge, and a domino is adjacent to the boundary
if it shares an edge with the boundary of the Aztec diamond. The 
\it north polar region \rm is defined to be the union of those north-going 
dominoes that are connected to the boundary by a sequence of adjacent 
north-going dominoes. The south, west and east polar regions are defined 
analogously. In this way a domino tiling partitions the Aztec diamond into 
four polar regions, where we have a regular brick wall pattern, and a fifth
central region, the \it temperate zone\rm, where the tiling pattern is 
irregular.

Consider the diagonal of white squares with opposite corners $Q_k^r$,
$k=0,\dots,n+1$, where $Q_k^r=(-r+k,n+1-k-r)$, $r=1,\dots,n$.
A \it zig-zag \rm path $Z_r$ in
$A_n$ from $Q_0^k$ to $Q_{n+1}^r$ is a path of edges going around these 
white squares. When going from $Q_k^r$ to $Q_{k+1}^r$ we can go either first
one step east and then one step south, or first one step south and then one 
step east. A domino tiling on $A_n$ defines a unique zig-zag path $Z_r$ from 
$Q_0^r$ to $Q_{n+1}^r$ if we require that the zig-zag path does not intersect
the dominoes. Similarly, we can define zig-zag paths from $P_0^r=
(-r,n-r)$ to $P_n^r=(n-r, -r)$ going around black squares.

We consider random tilings of the Aztec diamond, where each of the 
$2^{n(n+1)/2}$ possible tilings have the same probability. This
induces a probability
measure on the zig-zag paths. Consider a zig-zag path in
$A_n$ from $Q_0^k$ to $Q_{n+1}^r$ around white squares. Let $h_r<\dots
<h_1$ be those $k$ for which we go first east and then south when we go from
$Q_k^r$ to $Q_{k+1}^r$, $k=0,\dots, n$; there are exactly $r$ such $k$ if 
the zig-zag path comes from a domino tiling, [EKLP]. 
Call this zig-zag path $Z_r(h)$.
\nonumproclaim{Proposition 5.4} Let $\{h_1,\dots,h_r\}\subseteq\{0,\dots,n\}$ be the
positions of\break the 
east\/{\rm /}\/south turns in a zig-zag path $Z_r(h)$ in the Aztec
diamond $A_n$ from\break $(-r,n+1-r)$ to $(n+1-r,-r)$ around white squares{\rm .} Then{\rm ,} the
probability for this particular zig\/{\rm -}\/zag path is
$$
P[Z_r(h)]=\Bbb P_{\text{Kr},r,n,1/2}[h].\tag 5.12
$$
If $\{h_1,\dots,h_r\}\subseteq\{0,\dots,n-1\}$ are the positions of the 
south\/{\rm /}\/east turns in a zig\/{\rm -}\/zag path $Z'_r(h)$ in $A_n$ from
$(-r,n-r)$ to $(n-r,-r)$ around black squares{\rm ,} then
$$
P[Z'_r(h)]=\Bbb P_{\text{Kr},r,n-1,1/2}[h].\tag 5.13
$$
\endproclaim

\demo{Proof}
Let $\Cal U_r(h)$ be
the number of possible domino tilings above $Z_r(h)$ in the Aztec diamond.
From the arguments in [EKLP], see also [PS], it follows that
$$
\Cal U_r(h)=2^{r(r-1)/2}\prod_{1\le i<j\le r}\frac{h_i-h_j}{j-i}.\tag 5.14
$$
Let $k_1<\dots<k_{n+1-r}$ be defined by
$$
\{k_1,\dots,k_{n+1-r}\}=\{0,\dots,n\}\setminus\{h_1,\dots,h_r\}.
$$
If $\Cal L_r(h)$ is the number of domino tilings of the region below $Z_r(h)$
in $A_n$, then, using the symmetry of the Aztec diamond, we see that
$$
\Cal L_r(h)=2^{(n+1-r)(n-r)/2}\prod_{1\le i<j\le n+1-r}\frac{k_j-k_i}{j-i}.
\tag 5.15
$$
Thus, the probability for a certain zig-zag path $Z_r=Z_r(h)$, 
specified by $h$, is
$$
P[Z_r(h)]=\frac {2^{(n+1-r)(n-r)/2+r(r-1)/2}}{2^{n(n+1)/2}}
\prod_{1\le i<j\le r}\frac{h_i-h_j}{j-i}
\prod_{1\le i<j\le n+1-r}\frac{k_j-k_i}{j-i}.\tag 5.16
$$
If we let $h_i=\mu_i+r-i$, $1\le i\le r$ and $k_j=r+j-1-\mu_j'$, $1\le j\le
n+1-r$, then $\mu$ and $\mu'$ are conjugate partitions; compare with the proof of
Lemma 5.1 ($N=r$, $M=n+1-r$, $s_i=h_i+1$, $r_j=k_j+1$). We see from (5.16)
that
$$
P[Z_r(h)]=2^{-(n+1-r)r}\left(\prod_{j=1}^{r-1}\frac 1{j!}\right)
\left(\prod_{j=1}^{n-r}\frac 1{j!}\right)V_r(\mu)V_{n+1-r}(\mu'),
$$
where $\mu\in\Omega_r(n+1-r)$. We can now apply Lemma 5.1, which gives
$$
P[Z_r(h)]=\frac{2^{r(r-1)}}{(n!)^r}\prod_{j=1}^{r-1}\frac{(n-j)!}{j!}
V_r(\mu)^2\prod_{j=1}^r\binom{n}{\mu_j+r-j}\frac 1{2^n}.
$$
Now, $h_i=\mu_i+r-i$, so we obtain
$$
P[Z_r(h)]=\frac{2^{r(r-1)}}{(n!)^r}\prod_{j=1}^{r-1}\frac{(n-j)!}{j!}
\Delta_r(h)^2\prod_{j=1}^r\binom{n}{h_j}\frac 1{2^n},\tag 5.17
$$
which is the Krawtchouk ensemble. Note that in (5.17) the order of the $h_i$'s 
is unimportant, so we can let $\{h_1,\dots,h_r\}\subseteq\{0,\dots,n\}$ be 
the (unordered) positions of the east/south turns. A completely analogous 
argument applies to the zig-zag paths in $A_n$ from $P_0^r$ to $P_n^r$ around 
black squares. This completes the proof.
\enddemo

It is proved in [JPS] that, with probability tending to 1 as
$n\to\infty$, the asymptotic shape of the temperate zone is a circle
centered at the origin and tangent to the boundary of the Aztec
diamond ({\it the arctic circle} theorem). This can be deduced from
Proposition 5.4 and the general results in Section 4 of [Jo3]. The
arctic circle is determined by the endpoints of the support of the
equilibrium measure (or the points where it saturates). Also, from
Theorem 5.3, we see that the fluctuations of the temperate zone around
the arctic circle is described by the Tracy-Widom distribution. This
can also be deduced from the fact, derived in [JPS], that the shape of
a polar region is related to the shape of a randomly growing Young
diagram. The growth model obtained is exactly the discrete time growth
model studied in [Jo3], and we can apply the results of that
paper. See [Jo4] for more details.

Finally, we will shortly discuss another random tiling problem related
to plane partitions using the combinatorial analysis by Cohn, Larsen and
Propp in [CLP]. For more details and pictures see the paper [CLP]. 
Plane partitions in an $a\times b\times c$ box can be seen to be in 
one-to-one correspondence with tilings of an $a,b,c$-hexagon with unit
rhombi with angles $\pi/3$ and $2\pi/3$, called \it lozenges\rm. An
$a,b,c$-hexagon has sides of length $a,b,c,a,b,c$ (in clockwise order), equal angles and
the length of the horizontal sides is $b$. If the major diagonal of the 
lozenge is vertical we talk about a \it vertical lozenge\rm. Consider
the uniform distribution on the set of all possible tilings of the 
$a,b,c$-hexagon with lozenges, which corresponds to the uniform
distribution on all plane partitions in the $a\times b\times c$ box. 
For simplicity we will now restrict ourselves to the case $a=b=c$. A 
horizontal line $k$ steps from the top, $k=0,\dots,a$ will intersect the 
vertical lozenges at positions $h_1+1,\dots, h_k +1$,
$0\le h_1<\dots<h_k\le a+k-1$, otherwise it passes through sides of the 
lozenges. A random tiling induces a probability measure on the sequences
$h=(h_1,\dots,h_k)$. Interpreting the formulas in Theorem 2.2 in [CLP]
we see that the probability for $h$ is
$$
P[h]=\frac 1{Z_{k,a}}\Delta_k(h)^2\prod_{j=1}^k\binom{h_j+a-k}{h_j}
\binom{2a-1-h_j}{a+k-1-h_j},\tag 5.18
$$
where $Z_{k,a}$ is a constant that can be computed explicitly. Note that the 
measure is symmetric in the $h_i$'s so we can regard (5.18) as a measure on 
$\{0,\dots,a+k-1\}^k$. Thus, again we get a discrete orthogonal polynomial
ensemble, this time with the weight
$$
w(x)=\binom{x+\alpha}{x}\binom{N+\beta-x}{N-x}\tag 5.19
$$
on $\{0,\dots,N\}$, with $\alpha=\beta=a-k$ and $N=a+k-1$. The orthogonal
polynomials for this weight are the Hahn polynomials, [NSU], so (5.18) 
should be called the \it Hahn ensemble\rm. If we do not have $a=b=c$
we will again get a weight function of the form (5.19) but with
different values of $\alpha,\beta$ and $N$ and with a different number
of particles. 
This model is further discussed in [Jo4], but to obtain the
Tracy-Widom distribution in this model is more complicated due to the
fact that it is less straightforward to compute the asymptotics of the
Hahn polynomials.

\demo{Acknowledgement}
I thank Eric Rains, Craig Tracy
and Harold Widom for helpful conversations and correspondence. I also
thank Alexei Borodin, Andrei Okounkov and Grigori Olshanski for keeping
me informed about their work, and Timo Sepp\"al\"ainen for drawing my
attention to the papers [JPS] and [CEP]. Part of this work was done
while visiting MSRI and I would like to express my gratitude to its 
director David Eisenbud for inviting me and to Pavel Bleher and 
Alexander Its for organizing the program on Random Matrix Models and their
Applications. This work was supported by the Swedish Natural Science
Research Council (NFR).

\AuthorRefNames [BDJ1]
\references

[AD] \name{D. Aldous} and \name{P. Diaconis},  Longest increasing subsequences:
from patience sorting to the Baik-Deift-Johansson theorem,
{\it Bull.\ Amer.\ Math.\ Soc\/}.\ {\bf 36} (1999), 199--213.

[BR1] \name{J. Baik} and  \name{E. Rains},  Algebraic aspects of increasing
subsequences,\hfill\break  xxx.lanl.gov/abs/math.CO/9905083.

[BR2] \name{J. Baik}  and  \name{E. Rains},  The asymptotics of monotone
subsequences of involutions, xxx.lanl.gov/abs/math.CO/9905084.

[BDJ1] \name{J. Baik, P. A. Deift}, and \name{K. Johansson},  On the distribution
of the length of the longest 
increasing subsequence in a random permutations\rm, 
{\it J.\ Amer.\ Math.\ Soc\/}.\ {\bf 12} (1999), 1119--1178.

[BDJ2] \bibline,  On the distribution
of the length of the second row of a Young diagram under Plancherel
measure, {\it Geom.\ Func.\ Anal.\/} {\bf 10} (2000), 702--731.

[Bi] \name{P. Biane},  Representations of symmetric groups and free
probability, {\it Adv.\ Math\/}.\ {\bf 138} (1998), 126--181.

[BO1] \name{A. Borodin} and  \name{G. Olshanski},  Distributions on partitions, point
processes and the hypergeometric kernel, {\it Comm.\ Math.\ Phys\/}.\
{\bf 211} (2000), 335--338.

[BO2] \bibline,  $Z$-measures on partitions, 
Robinson-Schensted-Knuth correspondence, and $\beta=2$ random matrix
ensembles, xxx.lanl.gov/abs/math.CO/9905189.

[BoOk] \name{A. Borodin} and  \name{A. Okounkov},  A Fredholm determinant formula
for Toeplitz determinants, {\it Integral Equations Operator Theory\/}
{\bf 37} (2000), 386--396.

[BOO] \name{A. Borodin, A. Okounkov}, and \name{G. Olshanski}, Asymptotics of
Plancherel measures for symmetric groups, {\it J.\ Amer.\ Math.\ Soc\/}.\
{\bf 13} (2000), 481--515.

[CEP]  \name{H. Cohn, N. Elkies}, and \name{J. Propp},  Local statistics for random domino
tilings of the Aztec diamond, {\it Duke Math.\ J\/}.\ {\bf 85} 
(1996), 117--166.

[CLP] \name{H. Cohn, M. Larsen}, and \name{J. Propp}, The shape of a typical
boxed plane partition, {\it New York J.\ Math\/}.\ {\bf 4} (1998), 
137--165.

[DZ] \name{J.-D. Deuschel} and \name{O. Zeitouni}, On increasing subsequences of
I.I.D.\  samples, {\it Combin.\ Probab.\ Comput\/}.\ {\bf 8} (1999),
247--263.

[EKLP] \name{N. Elkies, G. Kuperberg, M. Larsen}, and \name{J. Propp}, 
Alternating-sign matrices and domino tilings (Part I),
{\it J.\ Algebraic Combin\/}.\ {\bf 1} (1992), 111--132.

[Fu] \name{W. Fulton}, {\it  Young Tableaux\/}, {\it London Math.\ Soc.\
Student Texts\/} {\bf 35}, Cambridge Univ.\ Press, Cambridge, 1997.

[Ge] \name{I. M. Gessel}, Symmetric functions and $P$-recursiveness,
{\it J.\ Combin.\ Theory\/} {\bf 53} (1990), 257--285.

[Go] \name{W. M. Y. Goh}, Plancherel-Rotach asymptotics for Charlier
polynomials, {\it Constr.\ Approx\/}.\ {\bf 14} (1998), 151--168.

[JPS] \name{W. Jockush, J. Propp}, and \name{P. Shor}, Random domino tilings
and the arctic circle theorem, preprint, 1995, xxx.lanl.gov/abs/math.CO/9801068.

[Jo1] \name{K. Johansson},  On fluctuations of Eigenvalues of random Hermitian
matrices,  {\it Duke Math.\ J\/}.\ {\bf 91} (1998), 151--204.

[Jo2] \bibline,  The longest increasing subsequence in a
random permutation and a unitary random matrix model,
{\it Math. Res. Lett\/}.\
{\bf 5} (1998), 63--82.

[Jo3] \bibline, Shape fluctuations and random matrices,
{\it Comm.\ Math.\ Phys\/}.\ {\bf 209} (2000), 437--476.

[Jo4] \bibline,  Non-intersecting paths, random tilings and
random matrices, xxx.lanl.gov\break /abs/math.PR/0011250.

[Ke1] \name{S. Kerov}, Transition probabilities of continual Young
diagrams and the Markov moment problem, {\it Funct.\ Anal.\ Appl\/}.\ 
{\bf 27} (1993), 104--117.

[Ke2] \bibline,  Asymptotics of the separation of roots of orthogonal
polynomials, 
 {\it St.\ Petersburg Math.\ J\/}.\ {\bf 5} (1994),
925--941.

[Kn] \name{D.\ E. Knuth},  Permutations, matrices and generalized Young
tableaux, {\it Pacific J.\ Math\/}.\ {\bf 34} (1970), 709--727.

[Me] \name{M.\ L.\ Mehta},  {\it Random Matrices}, 2nd ed., Academic Press, 
Boston, MA, 1991.

[NSU] \name{A. F. Nikiforov, S. K. Suslov}, and \name{V. B. Uvarov}, {\it Classical
Orthogonal Polynomials of a Discrete Variable}, {\it Springer Series in
Computat.\ Phys.\  Physics\/}, Springer-Verlag, New York, 1991.

[Ok] \name{A. Okounkov}, Random matrices and random permutations,\hfill\break  
xxx.lanl.gov/abs/math.CO/9903176.  

[PS] \name{J. Propp}  and \name{R. Stanley},  Domino tilings with barriers,
{\it J.\ Combin.\ Theory\/} Ser.\ A {\bf 87} (1999), 347--356.

[Ra] \name{E. Rains},  Increasing subsequences and the classical
groups, {\it Electron.\ J.\ Combin.\/} {\bf 5} (1998), 12pp.

[Re] \name{A. Regev},  Asymptotic values for degrees associated with
strips of Young diagrams, {\it Adv.\ in Math\/}.\ {\bf 41} (1981),
115--136.

[Sa] \name{B. Sagan}, \it The Symmetric Group\rm, Brooks/Cole 
Adv.\ Books and Software, Pacific Grove, CA,  1991.

[Se1] \name{T. Sepp\"al\"ainen},  Large deviations for increasing
subsequences on the plane, {\it Probab.\ Theory Related Fields\/} {\bf 112}
(1998), 221--244.

[Se2] \bibline,  Exact limiting shape for a simplified
model of first-passage percolation on the plane, {\it Ann.\ Probab\/}.\ 
{\bf 26}
(1998), 1232--1250.

[St] \name{R. P. Stanley}, {\it Enumerative Combinatorics,} Vol. 2,
{\it Cambridge Stud.\ in Adv.\ Math\/}.\ {\bf 62}, Cambridge Univ.\ Press, Cambridge, 1999.

[TW1] \name{C. A. Tracy} and \name{H. Widom},  Level spacing distributions and
the Airy kernel, {\it Comm.\ Math.\ Phys\/}.\ {\bf 159} (1994),
151--174. 

[TW2] \bibline,  Correlation functions, cluster
functions, and spacing distributions for random matrices,
{\it J.\ Statist.\ Phys\/}.\ {\bf 92} (1998), 809--835.

[TW3] \bibline,  On the distributions of the
lengths of the longest monotone subsequences in random words,
xxx.lanl.gov/abs/math.CO/9904042.

[Ve] \name{A. Vershik},  Asymptotic combinatorics and algebraic
analysis, {\it Proc.\ Internat.\ Congress of Mathematicians\/}, 
(Z\"urich, 1994), 1384--1394, Birkh\"auser Basel, 1995.

[VK] \name{A. Vershik} and \name{S. Kerov}, Asymptotics of the Plancherel
measure of the symmetric group and the limiting form of Young
tables, {\it Soviet Math.\ Dokl\/}.\ {\bf 18} (1977), 527--531.

[Wa] G. N. Watson\name{}, {\it A Treatise on the Theory of Bessel
Functions\/}, Cambridge Univ.\ Press, Cambridge, 1952.

\endreferences
\bye bye!

%\flushpar \heading REFERENCES\endheading \medskip

%[AD} D. Aldous, P. Diaconis, {\it Longest increasing subsequences:
%From patience sorting to the Baik-Deift-Johansson theorem,}
%Bull. AMS, {\bf 36} (1999), 199 - 213

%[BR1} J. Baik, E. Rains, \it Algebraic aspects of increasing
%subsequences\rm,\newline  xxx.lanl.gov/abs/math.CO/9905083

[BR2} J. Baik, E. Rains, \it The asymptotics of monotone
subsequences of involutions\rm, xxx.lanl.gov/abs/math.CO/9905084

[BDJ1} J. Baik, P. A. Deift, K. Johansson, \it On the distribution
of the length of the longest 
increasing subsequence in a random permutation\rm, 
J. Amer. Math. Soc., {\bf 12}, (1999), 1119 - 1178

[BDJ2} J. Baik, P. A. Deift and K. Johansson, \it On the distribution
of the length of the second row of a Young diagram under Plancherel
measure\rm, to appear in Geometric and Funct. Anal.

[Bi} P. Biane, \it Representations of symmetric groups and free
probability\rm, Adv. Math., {\bf 138}, (1998), 126 - 181

[BO1} A. Borodin, G. Olshanski, \it Distributions on Partitions, Point
Processes and the Hypergeometric Kernel\rm, xxx.lanl.gov/abs/math.CO/9904010

[BO2} A. Borodin, G. Olshanski, \it Z-measures on partitions, 
Robinson-Schensted-Knuth correspondence, and $\beta=2$ random matrix
ensembles\rm, \newline xxx.lanl.gov/abs/math.CO/9905189

[BoOk} A. Borodin, A. Okounkov, \it A Fredholm determinant formula
for Toeplitz determinants\rm, math.CA/9907165

[BOO} A. Borodin, A. Okounkov, G. Olshanski, \it On asymptotics of
Plancherel measures for symmetric groups\rm, xxx.lanl.gov/abs/math.CO/990532

[CEP} H. Cohn, N. Elkies, J. Propp, \it Local statistics for random domino
tilings of the Aztec diamond\rm, Duke Math. J., {\bf 85}, (1996), 117 - 166

[CLP} H. Cohn, M. Larsen, J. Propp, \it The shape of a typical
boxed plane partition\rm, New York J. of Math., {\bf 4}, (1998), 137 -
165

[DZ} J.-D. Deuschel, O. Zeitouni, \it On increasing subsequences of
i.i.d. samples\rm, Combin. Probab. Comput., {\bf 8}, (1999), 247 - 263

[EKLP} N. Elkies, G. Kuperberg, M. Larsen, J. Propp, \it
Alternating-Sign Matrices and Domino Tilings (Part I)\rm,
J. of Algebraic Combin., {\bf 1}, (1992), 111- 132

[Fu} W. Fulton, \it Young Tableaux\rm, London Mathematical Society,
Student Texts 35, Cambridge Univ. Press, 1997

[Ge} I. M. Gessel, \it Symmetric functions and P-recursiveness\rm,
J. Combin. Theory Ser. A, {\bf 53}, (1990), 257 - 285

[Go} W. M. Y. Goh, \it Plancherel-Rotach asymptotics for Charlier
polynomials\rm, Constr. Approx., {\bf 14}, (1998), 151 - 168

[JPS} W. Jockush, J. Propp, P. Shor, \it Random domino tilings
and the arctic circle theorem\rm, preprint 1995, xxx.lanl.gov/abs/math.CO/9801068

[Jo1} K. Johansson, \it On Fluctuations of Eigenvalues of Random Hermitian
Matrices\rm,  Duke Math. J., {\bf 91}, (1998), 151 - 204

[Jo2} K. Johansson, \it The longest increasing subsequence in a
random permutation and a unitary random matrix model\rm,
Math. Res. Lett.,
{\bf 5}, (1998), 63 - 82

[Jo3} K. Johansson, \it Shape fluctuations and random matrices\rm,
xxx.lanl.gov/abs/math.CO/9903134, Commun. Math. Phys., {\bf 209}, (2000), 437 - 476

[Jo4} K. Johansson, \it Non-intersecting paths, random tilings and
random matrices\rm, in preparation

[Ke1} S. Kerov, \it Transition probabilities of continual Young
diagrams and the Markov moment problem\rm, Func. Anal. Appl., {\bf
27}, (1993), 104 - 117

[Ke2} S. Kerov, \it The asymptotics of interlacing roots of
orthogonal polynomials\rm, St. Petersburg Math. J., {\bf 5}, (1994),
925 - 941

[Kn} D. E. Knuth, \it Permutations, Matrices and Generalized Young
Tableaux\rm, Pacific J. Math., {\bf 34}, (1970), 709 - 727

[Me} Mehta, M. L., {\it Random Matrices,} 2nd ed., Academic Press, 
San Diego 1991

[NSU} A. F. Nikiforov, S. K. Suslov, V. B. Uvarov, \it Classical
Orthogonal Polynomials of a Discrete Variable\rm, Springer Series in
Computational Physics, Springer-Verlag, Berlin Heidelberg, 1991

[Ok} A. Okounkov, \it Random Matrices and Random Permutations\rm,
\newline xxx.lanl.gov/abs/math.CO/9903176  

[PS} J. Propp, R. Stanley, \it Domino tilings with barriers\rm,
(1998), J. Combin. Th. Ser.A, {\bf 87}, 347 - 356

[Ra} E. Rains, \it Increasing subsequences and the classical
groups\rm, Electr. J. Of Combinatorics\rm, {\bf 5(1)}, (1998), R12

[Re} A. Regev, \it Asymptotic values for degrees associated with
strips of Young diagrams\rm, Adv. in Math., {\bf 41}, (1981), 115 - 136

[Sa} B. Sagan, \it The Symmetric Group\rm, Brooks/Cole Publ. Comp., 1991

[Se1} T. Sepp\"al\"ainen, \it Large Deviations for Increasing
Subsequences on the Plane\rm, Probab. Theory Relat. Fields, {\bf 112},
(1998), 221 - 244

[Se2} T. Sepp\"al\"ainen, \it Exact limiting shape for a simplified
model of first-passage percolation on the plane\rm, Ann. Prob., {\bf 26},
1232 - 1250

[St} R. P. Stanley, {\it Enumerative Combinatorics,} Vol. 2,
Cambridge University Press, 1999

[TW1} C. A. Tracy, H. Widom, \it Level Spacing Distributions and
the Airy Kernel\rm, Commun. Math. Phys., {\bf 159}, (1994), 151 - 174 

[TW2} C. A. Tracy, H. Widom, \it Correlation Functions, Cluster
Functions, and Spacing Distributions for Random Matrices\rm,
J. Statist. Phys., {\bf 92}, (1998), 809 - 835

[TW3} C. A. Tracy, H. Widom, \it On the distributions of the
lengths of the longest monotone subsequences in random words\rm,
xxx.lanl.gov/abs/math.CO/9904042

[Ve} A. Vershik, \it Asymptotic Combinatorics and Algebraic
Analysis\rm, Proceedings of the International Congress of Mathematics,
Z\"urich, Switzerland  1994, \linebreak 
Birkh\"auser Verlag, Basel, 1995, 1384 -
1394.

[VK} A. Vershik, S. Kerov, \it Asymptotics of the Plancherel
measure of the symmetric group and the limiting form of Young
tables\rm, Soviet Math. Dokl., {\bf 18}, (1977), 527-531

[Wa} G. N. Watson, \it A treatise on the theory of Bessel
functions\rm, Cambridge University Press, Cambridge, 1952

\end